\documentclass[sensors,article,submit,pdftex,moreauthors]{Definitions/mdpi} 
\usepackage{lmodern}
\usepackage{anyfontsize}
\firstpage{1} 
\pubvolume{1}
\issuenum{1}
\articlenumber{0}
\pubyear{2023}
\copyrightyear{2023}
\datereceived{ 14 March 2023 } 
\daterevised{13 April 2023 } 
\dateaccepted{14 April 2023 } 
\datepublished{ 17 April 2023} 
\hreflink{https://doi.org/} 

\usepackage{amsmath,amssymb,amsfonts,latexsym,graphicx}
\usepackage{bm}

\usepackage{algorithm}
\usepackage{algorithmic}

\usepackage{graphicx}
\usepackage{subcaption}
\usepackage{changes}

\graphicspath{ {images/} }

\DeclareMathOperator{\Cov}{Cov}
\DeclareMathOperator{\diag}{diag}

\Title{An optimal linear fusion estimation algorithm of reduced dimension for  \texorpdfstring{$\mathbb{T}$}-proper systems with multiple packet dropouts}

\TitleCitation{An optimal linear fusion estimation algorithm of reduced dimension for T-proper systems with multiple packet dropouts.}

\Author{Rosa M.     Fern\'andez-Alcal\'a $^{*}$ \orcidA{}, José D.     Jiménez-López \orcidB{}, Nicolas Le Bihan \orcidC{}, and  Clive Cheong Took \orcidC{}}

\AuthorNames{Rosa M.     Fern\'andez-Alcal\'a, José D.     Jiménez-López, Nicolas Le-Bihan, and Clive Cheong Took}

\address{%
$^{1}$ \quad Department of Statistics and Operations Research, University of Jaén, Paraje Las Lagunillas, 23071, Jaén, Spain;
 rmfernan@ujaen.    es (R.    M.    F.    -A.    ); jdomingo@ujaen.    es (J.    D.    J.    -L.    )  \\
 
$^{2}$ \quad Department of  Images and Signals, CNRS/GIPSA-Lab, 38402 Saint Martin d’He`res Cedex, France;  nicolas.le-bihan@gipsa-lab.grenoble-inp.  fr (N. L. B. )  

$^{3}$ \quad Department of Electronic Engineering. Royal Holloway, London, TW20 OEX, UK; Clive.CheongTook@rhul.ac. uk (C.    C.    T.    )\\

}

\corres{Correspondence: rmfernan@ujaen.    es (R.    M.    F.    -A.    ); Tel.    : +34-953-212-449}

\abstract{This paper analyses the centralized fusion linear estimation problem in multi-sensor systems with multiple packet dropouts and correlated noises. Packet dropouts are modeled by independent Bernoulli distributed random variables.  This problem is addressed in the tessarine domain under conditions of \texorpdfstring{$\mathbb{T}_1$} and \texorpdfstring{$\mathbb{T}_2$}-properness, which entails a reduction in the dimension of the problem and,  consequently, computational savings. The methodology proposed enables us to provide an optimal (in the least-mean-squares sense) linear fusion filtering algorithm for estimating the tessarine state with a lower computational cost than the conventional one devised in the real field.  Simulation results illustrate the performance and advantages of the solution proposed in different settings.}


\begin{document}


\section{Introduction}

In sensor networks, the problem of estimating the state observed by multiple sensors has been analyzed extensively in recent decades due to the variety of applications they have in signal processing (see, e.g., \cite{Kurkin2017, Kondru2019, Huang2020, Liping2021, Tian2022, Kong2023, Zewge2023, Zheng2023, Liu2023}). 

In networked systems, sensor failures, network congestion, communications interferences or noises can cause random packet dropouts in data transmissions and, consequently, it is possible that the measurements available for state estimation are not always updated. These packet dropouts can be described by stochastic parameter systems that define the strategy followed to compensate for packet loss \cite{Schenato2007, Sun2008, Sun2014, Zhao2016, Lu2019, Song2022, Li2022a, Su2022,Geng2022, Xu2023, Zhu2023, Ni2023, Elahi2023}.  

For multi-sensor systems, the potential of fusion estimation techniques to produce consistent and accurate estimators has been demonstrated. Thus, these techniques have also been applied to multi-sensor systems with multiple packet dropouts, giving rise to centralized as well as distributed fusion estimation algorithms (see, e.g., \cite{Schenato2007, Ma2011, raguila2012, 
Ma2013, Li2014, Liping2021, Li2022}). In general, centralized fusion methodology yields optimal estimators, but the computational load involved can be a handicap in practical applications.

Alternatively, 4D hypercomplex-based signal processing has been satisfactorily applied as a dimension  reduction approach in multi-sensor fusion estimation problems with uncertainties \cite{Yuan2015, Talebi2016, Tannous2016, Navarro2019, Wu2019, Talebi2020, Fernandez2021, Jimenez2021}. Effectively, the benefit of using hypercomplex algebras is twofold: firstly, they may provide a compact representation of multidimensional signals and a better insight into the structure of the problem than that provided by a traditional or real formalism, and secondly, the characterization of certain properness properties related to the vanishing of some correlation or pseudo correlation functions, makes the dimension of the processes involved may be reduced. Then, even though the optimal processing in the 4D hypercomplex field is the widely linear (WL) processing, which implies to work on a four-dimensional  vector; under properness conditions,  WL processing is equivalent to using a signal processing based on a vector of reduced dimension. Notice that there is not an algebra that always leads to the better solution, but the choice of the best algebra at each situation depends on the proper characteristics of the processes involved.

From among the different 4D hypercomplex structures, quaternions and, more recently, tessarines have been the most popular algebras used in signal processing. They are characterized by different multiplication rules that endorse them different algebraic properties of interest: on the one hand, quaternions are a noncommutative division algebra and, on the other hand,  tessarines are a commutative non-division algebra.  Nevertheless,  the fact that they have zero divisors does not have a major effect on their practical applications. Recently, due to the advantage of working with commutative algebras, the multi-sensor fusion estimation problem for systems with uncertain measurements has been addressed for tessarine signals with properness properties. Specifically, under $\mathbb{T}_1$ and $\mathbb{T}_2$ properness conditions, Kalman filter-like centralized and distributed fusion estimation algorithms have been proposed in \cite{Fernandez2021, Jimenez2021} by considering different uncertainty situations (missing measurements and/or random delays), as well as correlated noises. The main interest of these algorithms lies in the reduction of the computational burden they entail under properness conditions, as achieving 
this computational saving from a real formalism is not possible. Nevertheless, the benefits of this methodology still have not been exploited in multiple sensor stochastic systems with packet dropouts.  

This paper deals with the linear least-mean-squares (LLMS) fusion filtering problem for multi-sensor $\mathbb{T}_k$-proper, ($k=1,2$)  tessarine systems with multiple packet dropouts. At each sensor, the multiple packet dropouts are described by independent Bernoulli distributed tessarine random vectors that, at any instant of time $t$, indicates whether the measurement output is received or lost, and in this second situation, the latest measurement available is used. Moreover,  our formulation of the problem includes a possible correlation between the state and measurement noises. In this setting, using centralized fusion Kalman filter technique, and based on a  $\mathbb{T}_k$-proper signal processing,  an optimal linear fusion filtering algorithm of reduced dimension is provided for estimating the state as well as its mean squared error. Additionally, the performance of the solution proposed and its superiority over its counterpart in the quaternion domain is experimentally analyzed for the cases of $\mathbb{T}_1$-properness as well as $\mathbb{T}_2$-properness,   by using a numerical example. In summary, the challenges of this paper are described in brief by the following items: 1) to address the LLMS fusion filtering problem  for multisensor systems with multiple packet dropouts and correlated noises in the tessarine domain, 2) to establish conditions on the state space sytem that guarantee the $\mathbb{T}_k$-properness of the processes involved, 3) to analyze the implications of $\mathbb{T}_k$-properness in the reduction of the dimension of the problem, 4) to derive a recursive algorithm to obtain the optimal fusion filters, 5) to illustrate numerically the benefits of the proposed solution over their counterparts in the quaternion setting, under $\mathbb{T}_k$-properness conditions.

The paper is organized as follows: In Section 2  the main concepts and properties in the tessarine domain are reviewed. Specifically, a tessarine random signal vector is defined, its conjugations, the real and augmented vectors, the {\em pseudo} auto and cross correlation functions, the $\mathbb{T}_k$-properness, and the $\star$ product between tessarines. Section 3 describes the centralized fusion filtering problem for multi-sensor systems with multiple packet dropouts and its formulation in the tessarine domain under conditions of $\mathbb{T}_k$-properness. On the basis of the state-space model of the tessarine signal and the observations, the WL stacked state-space system is built from the augmented vectors of the processes involved and, afterwards, under $\mathbb{T}_k$-properness conditions, an equivalent form of reduced dimension for the available observation equation is presented. Next, the $\mathbb{T}_k$-proper centralized linear fusion filtering algorithm based on Kalman filter techniques  is presented in Section 4. Note that, to preserve the continuity of exposition, the derivation of the formulas of this algorithm has been deferred to Appendix A. In Section 5, two numerical examples, one of them over simulated values, and the other one on a realistic phenomena,  illustrate the theoretical results obtained. The paper finishes with the concluding remarks in Section 6. 

\subsection{Notations}

The following standard notation is used throughout this paper:  scalars are denoted by lightface letters, while boldface lowercase and boldface uppercase letters represent the vectors and matrices, respectively.  The symbol  $\mathbf{0}_{n\times m}$ (respectively, $\mathbf{0}_n$)  stands for the $n\times m$ matrix (respectively, $n$ column vector) whose all elements are zeros, $\mathbf{I}_n$ is the $n\times n$ identity matrix, and  $\mathbf{1}_n$  denotes the $n$ column vector of ones.    

$\mathbb{Z}$, $\mathbb{R}$ and $\mathbb{T}$ denote, respectively, the set of integer, real, and tessarine numbers. $\mathbb{R}^{n}$ (respectively, $\mathbb{T}^{n}$) is the set of all $n$-dimensional real (respectively, tessarine) vectors, and $\mathbb{R}^{n\times m}$ (respectively, $\mathbb{T}^{n\times m}$) refers to the set of all $n\times m$-dimensional real (respectively, tessarine) matrices. Moreover, the superscripts ``*'', ``$\texttt{T}$'', and ``$\texttt{H}$'' symbolize the tessarine conjugate, transpose, and Hermitian transpose, respectively.

The notation $E[\cdot]$ represents the mathematical expectation,  ${\Cov}(\cdot)$ is the covariance operator, and ${\diag}(\cdot)$ denotes the diagonal (or block diagonal) matrix with the input arguments on the main diagonal.  Finally, $\delta_{t,s}$ represents the Kronecker delta function, and the Hadamard and Kronecker product operators are symbolized by ``$\circ$'' and ``$\otimes$'', respectively.

\section{Definitions and Preliminaries}

This section is devoted to stating the core concepts and results in the tessarine domain that will be used throughout the paper.

Unless otherwise indicated, we shall assume that all random variables have zero mean.

\begin{Definition}
A tessarine random signal vector $\mathbf{x}(t)\in \mathbb{T}^n$ is a four-dimensional hypercomplex stochastic process defined as \cite{Navarro2020}
  \begin{equation*}\label{tessarine random signal}
    \mathbf{x}(t)=\mathbf{x}_r(t)+\eta\mathbf{x}_{\eta}(t)+\eta'\mathbf{x}_{\eta'}(t)+\eta''\mathbf{x}_{\eta''}(t),  \ \ t\in \mathbb{Z},
  \end{equation*}
with $\mathbf{x}_{\nu}(t)\in\mathbb{R}^n$, for $\nu=r,\eta,\eta',\eta''$, and where $\{\eta,\eta',\eta''\}$ are  hyper-imaginary units such that:
  \begin{equation*}
    \eta\eta'=\eta'',\quad \eta'\eta''=\eta,\quad \eta''\eta=-\eta',\quad \eta^2=\eta''{\,^2}=-1,\quad \eta'{\,^2}=1.    
  \end{equation*}
\end{Definition}

Consider $\mathbf{x}(t), \mathbf{y}(t) \in \mathbb{T}^n$, tessarine random signal vectors given in Definition 1. The following concepts and properties can be established.   
Let $\boldsymbol{\Gamma}_{\mathbf{x}}(t,s)=E[\mathbf{x}(t) \mathbf{x}^{\texttt{H}}(s)]$  be the {\it pseudo} autocorrelation function of $\mathbf{x}(t)\in \mathbb{T}^n$, and  $\boldsymbol{\Gamma}_{\mathbf{x}\mathbf{y}}(t,s)=E[\mathbf{x}(t)
\mathbf{y}^{\texttt{H}}(s)]$ the {\it pseudo} cross-correlation function of  $\mathbf{x}(t), \mathbf{y}(t) \in \mathbb{T}^n$, $\forall t, s \in \mathbb{Z}$.

In the tessarine domain, the second-order statistical properties of $\mathbf{x}(t)\in \mathbb{T}^n$, are completely described from the augmented tessarine
signal vector  
\begin{equation}\label{augmented}
\bar{\mathbf{x}}(t)=[\mathbf{x}^{^\texttt{T}}(t), \mathbf{x}^{*^{\texttt{T}}}(t), \mathbf{x}^{\eta^{\texttt{T}}}(t),\mathbf{x}^{\eta''^{\texttt{T}}}(t)]^{\texttt{T}},
\end{equation}
where $\mathbf{x}^{*}(t)$ is the conjugate of $\mathbf{x}(t)$ defined as 
$$\mathbf{x}^{*}(t)=\mathbf{x}_{r}(t)-\eta \mathbf{x}_{\eta}(t)+\eta' \mathbf{x}_{\eta'}(t)-\eta''\mathbf{x}_{\eta''}(t),$$ and 
$$\mathbf{x}^{\eta}(t)=\mathbf{x}_{r}(t)+\eta \mathbf{x}_{\eta}(t)-\eta' \mathbf{x}_{\eta'}(t)-\eta''\mathbf{x}_{\eta''}(t),$$
$$\mathbf{x}^{\eta''}(t)=\mathbf{x}_{r}(t)-\eta \mathbf{x}_{\eta}(t)-\eta' \mathbf{x}_{\eta'}(t)+\eta''\mathbf{x}_{\eta''}(t).$$

Let $\mathbf{x}^r(t)=[\mathbf{x}_{r}^{\texttt{T}}(t),\mathbf{x}_{\eta}^{\texttt{T}}(t),\mathbf{x}_{\eta'}^{\texttt{T}}(t),\mathbf{x}_{\eta''}^{\texttt{T}}(t)]^\texttt{T}$ be  the real vector formed by the components $\mathbf{x}_{\nu}(t)\in\mathbb{R}^n$, $\nu=r,\eta,\eta',\eta''$, of $\mathbf{x}(t)\in \mathbb{T}^n$. The following relationship can be established:
\begin{equation*}\label{q2r}
\bar{\mathbf{x}}(t)=2\boldsymbol{\mathcal{T}}\mathbf{x}^r(t),
\end{equation*}
where $\boldsymbol{\mathcal{T}}=\frac{1}{2} \boldsymbol{\mathcal{A}}\otimes \boldsymbol{I}_n$, with
$$\boldsymbol{\mathcal{A}}= \left[
\begin{array}{rrrr}
1 &  \eta &  \eta'   &   \eta'' \\
1 &  -\eta &   \eta'  & -\eta''  \\
1 &  \eta &  -\eta'   &  -\eta''  \\
1 &  -\eta &  -\eta'   &   \eta'' \\
\end{array}\right].$$
Notice that  $\boldsymbol{\mathcal{T}}^{\texttt{H}}\boldsymbol{\mathcal{T}}=\boldsymbol{I}_{4n}$.


It should be highlighted that the properness profile of a tessarine random signal plays a key role in the choice of the suitable type of linear processing that leads to a reduction in the dimension of the problem. This properness profile is characterized by the degree of correlation between the imaginary components and the real component. In particular, two interesting types of properness can be defined in the tessarine domain \citep{Navarro2020, Navarro2021}.

\begin{Definition} 
Let $\mathbf{x}(t)\in \mathbb{T}^n$ be a tessarine random signal vector. It is said that:
\begin{itemize}
\item $\mathbf{x}(t)$ is $\mathbb{T}_1$-proper if, and only if, $\boldsymbol{\Gamma}_{\mathbf{x}\mathbf{x}^{\nu}}(t,s)$=0, for $\nu=*,\eta,\eta''$, and $\forall t, s \in \mathbb{Z}$,
\item $\mathbf{x}(t)$ is $\mathbb{T}_2$-proper if, and only if, $\boldsymbol{\Gamma}_{\mathbf{x}\mathbf{x}^{\nu}}(t,s)$=0, for $\nu=\eta,\eta''$, and $\forall t, s \in \mathbb{Z}$.
\end{itemize} 
Likewise, let $\mathbf{x}(t)\in \mathbb{T}^{n_1}$ and $\mathbf{y}(t)\in \mathbb{T}^{n_2}$ be two tessarine random signal vectors. It is said that:
\begin{itemize}
\item $\mathbf{x}(t)$ and $\mathbf{y}(t)$ are cross  $\mathbb{T}_1$-proper,  if, and only if, $\boldsymbol{\Gamma}_{\mathbf{x}\mathbf{y}^{\nu}}(t,s)$=0, for $\nu=*,\eta,\eta''$ , and $\forall t, s \in \mathbb{Z}$,
\item $\mathbf{x}(t)$ and $\mathbf{y}(t)$ are cross  $\mathbb{T}_2$-proper,  if, and only if, $\boldsymbol{\Gamma}_{\mathbf{x}\mathbf{y}^{\nu}}(t,s)$=0, for $\nu=\eta,\eta''$ , and $\forall t, s \in \mathbb{Z}$,
\item $\mathbf{x}(t)$ and $\mathbf{y}(t)$ are jointly $\mathbb{T}_k$-proper, for $k=1,2$, if, and only if, they are $\mathbb{T}_k$-proper and cross $\mathbb{T}_k$-proper.

\end{itemize}

\end{Definition} 

\begin{Remark}
In the tessarine domain, the optimal linear processing, the widely linear (WL) processing, is based on augmented tessarine vector of dimension $4n$ of the form given in \eqref{augmented}.   Nevertheless, when $\mathbb{T}_k$-properness conditions are satisfied,  the WL estimators coincides with the one obtained from a $\mathbb{T}_k$-proper linear processing, which uses only the information provided by the processes involved (case $k=1$) or the $2n$-dimensional augmented vector formed by the signal and its conjugate (case $k=2$). In consequence, $\mathbb{T}_k$-properness means there is a significant reduction in the dimension of the processes involved  \citep{Navarro2021}. 

\end{Remark}

Finally, a new product between two tessarine signal vectors is defined.
\begin{Definition}
Consider $\mathbf{x}(t),\mathbf{y}(s)\in\mathbb{T}^n$. The product $\star$ is defined by the expression
\begin{equation*}\label{star product}
 \mathbf{x}(t) \star  \mathbf{y}(s)=\mathbf{x}_r(t) \circ \mathbf{y}_r(s)+\eta \mathbf{x}_{\eta}(t)\circ \mathbf{y}_{\eta}(s)+\eta' \mathbf{x}_{\eta'}(t)\circ \mathbf{y}_{\eta'}(s)+\eta'' \mathbf{x}_{\eta''}(t)\circ \mathbf{y}_{\eta''}(s).       
 \end{equation*} 
\end{Definition}

Note that,  given  two random tessarine signal vectors  $\mathbf{x}(t),\mathbf{y}(s)\in\mathbb{T}^n$,  the augmented vector of $\mathbf{x}(t)\star\mathbf{y}(s)$ is $\overline{\mathbf{x}(t)\star\mathbf{y}(s)}=\boldsymbol{\mathcal{D}}^{{\mathbf{x}}}(t)\bar{\mathbf{y}}(s)$, with $\boldsymbol{\mathcal{D}}^{{\mathbf{x}}}(t)=\boldsymbol{\mathcal{T}}\diag(\mathbf{x}^r(t))\boldsymbol{\mathcal{T}}^{\texttt{H}}$.    

\section{Problem formulation}
Let $\mathbf{x}(t)\in \mathbb{T}^n$ be an n-dimensional tessarine state vector which is assumed to be observed from $R$ sensors perturbed by different additive noises  according to the state-space model:
\begin{equation*}\label{SS}
\begin{split}
\mathbf{x}(t+1)=&  \mathbf{F}_1(t)\mathbf{x}(t)+\mathbf{F}_2(t)\mathbf{x}^{*}(t)+\mathbf{F}_3(t)\mathbf{x}^{\eta}(t)+
\mathbf{F}_4(t)\mathbf{x}^{\eta''}(t)+\mathbf{u}(t), \ \ t\geq 0, \\
\mathbf{z}^{(i)}(t)=&  \mathbf{x}(t)+\mathbf{v}^{(i)}(t), \ \ t\geq 1, \ \ i=1,\ldots,R, 
\end{split}
\end{equation*}
with 
\begin{itemize}
\item  $\mathbf{F}_j(t)\in \mathbb{T}^{n\times n}$, $j=1,\ldots,4$: deterministic tessarine matrices.
\item $\mathbf{u}(t) \in \mathbb{T}^{n}$: tessarine white noises with {\it pseudo} variances $\mathbf{Q}(t)$.
\item $\mathbf{v}^{(i)}(t) \in \mathbb{T}^{n}$: tessarine white noises with  {\it pseudo} variances $\mathbf{R}^{(i)}(t)$.
\item $\mathbf{u}(t)$,$\mathbf{v}^{(i)}(t)$ correlated with $\boldsymbol{\Gamma}_{\mathbf{u}\mathbf{v}^{(i)}}(t,s)=\mathbf{S}^{(i)}(t) \delta_{t,s}$. 
 \item $\mathbf{v}^{(i)}(t)$, $\mathbf{v}^{(j)}(t)$ independent for any two sensors $i\neq j$.
 \item $\mathbf{x}(0)$ uncorrelated with  $\mathbf{u}(t)$ and $\mathbf{v}^{(i)}(t)$, for $ t\geq 0$, $i=1,\dots,R$.
 \item $\boldsymbol{\Gamma}_{\mathbf{x}}(0,0)=\mathbf{P}_0$. 
\end{itemize}

The packets or measured outputs $\mathbf{z}^{(i)}(t)$ are assumed to be affected by random packet dropouts characterized by Bernoulli distributed random variables that can be described by the following model: 
\begin{equation}\label{available observation equation}
 \begin{aligned}
   \mathbf{y}^{(i)}(t)&=  \boldsymbol{\gamma}^{(i)}(t)\ \star\ \mathbf{z}^{(i)}(t)
    +(\mathbf{1}_n-\boldsymbol{\gamma}^{(i)}(t))\star\mathbf{y}^{(i)}(t-1),\quad t\geq 2;
  \end{aligned}
 \end{equation}
 for $ i=1,\dots,R$, with  $\mathbf{y}^{(i)}(1)=  \mathbf{z}^{(i)}(1)$, and the $\star$ product given in Definition 3.   Moreover, at each sensor $ i=1,\dots,R$, the tessarine random vector  $ \boldsymbol{\gamma}^{(i)}(t)=[\gamma_{1}^{(i)}(t), \dots,\gamma_{n}^{(i)}(t) ]^{\texttt{T}}\in \mathbb{T}^{n}$ is of the form $\gamma_{j}^{(i)}(t)=\gamma_{j,r}^{(i)}(t)+\eta \gamma_{j,\eta}^{(i)}(t)+\eta' \gamma_{j,\eta'}^{(i)}(t)+\eta'' \gamma_{j,\eta''}^{(i)}(t)$, for $j=1,\ldots, n$, where  $\gamma_{j,\nu}^{(i)}(t)$ are  independent Bernoulli random variables with known probabilities $p_{j,\nu}^{(i)}(t)$, for $j=1,\dots,n$ and $\nu= r, \eta, \eta', \eta''$,  that indicates whether the 
corresponding component of the packet or measured output $\mathbf{z}^{(i)}(t)$ of sensor $i$ is received at time $t$ ($\gamma_{j,\nu}^{(i)}(t)=1$) or it is lost and the latest received previously component, corresponding to the measured output $\mathbf{z}^{(i)}(t-1)$  
is used at time $t$ ($\gamma_{j,\nu}^{(i)}(t)=0$).
Additionally, $\boldsymbol{\gamma}^{(i)}(t)$ and $ \boldsymbol{\gamma}^{(i)}(s)$  are assumed to be independent for $t\neq s$, and $ \boldsymbol{\gamma}^{(i)}(t)$ is independent of $\mathbf{x}(t)$, $\mathbf{u}(t)$, $\mathbf{v}^{(l)}(t)$ and $\boldsymbol{\gamma}^{(l)}(t)$, for $i\neq l$,  with $i,l=1,\dots,R$. 

\begin{Remark}
Observe that model \eqref{available observation equation} always considers the latest measurement output received when the current measurement output is lost during transmission. Hence, this model can be used to describe multiple packet dropouts.
\end{Remark}

\begin{Remark}\label{remark00}
Under the hypothesis established for the Bernoulli random variables $\gamma_{j,\nu}^{(i)}(t)$, it is not difficult to check that
\begin{align*}
E\left[{\gamma}^{(i_1)}_{j_1, \nu_1}(t){\gamma}^{(i_2)}_{j_2, \nu_2}(t)\right]&=\left\{
\begin{array}{l}
{p}^{(i_1)}_{j_1, \nu_1}(t),\ \  if\  i_1=i_2, j_1=j_2, \nu_1=\nu_2\\
{p}^{(i_1)}_{j_1, \nu_1}(t){p}^{(i_2)}_{j_2, \nu_2}(t), \ \  otherwise,
\end{array}
\right. 
\\
E\left[\left( 1-{\gamma}^{(i_1)}_{j_1, \eta_1}(t)\right) \left( 1-{\gamma}^{(i_2)}_{j_2, \eta_2}(t)\right) \right]& =\left\{
\begin{array}{l}
1-{p}^{(i_1)}_{j_1, \nu_1}(t), \ \  if\  i_1=i_2, j_1=j_2, \nu_1=\nu_2\\
\left(1-{p}^{(i_1)}_{j_1, \nu_1}(t)\right)\left(1-{p}^{(i_2)}_{j_2, \nu_2}(t)\right), \ \  otherwise,
\end{array}
\right. 
\\
E\left[{\gamma}^{(i_1)}_{j_1, \eta_1}(t) \left( 1-{\gamma}^{(i_2)}_{j_2, \eta_2}(t)\right) \right]& =\left\{
\begin{array}{l}
0, \ \  if\  i_1=i_2, j_1=j_2, \nu_1=\nu_2\\
{p}^{(i_1)}_{j_1, \nu_1}(t) \left(1-{p}^{(i_2)}_{j_2, \nu_2}(t)\right), \ \  otherwise,
\end{array}
\right. 
\end{align*}
for any $j_1,j_2=1,\ldots, n$, $\nu_1,\nu_2=r, \eta,\eta',\eta''$ and $i_1,i_2=1,\ldots, R$.
\end{Remark}

In this setting, and based on the information supplied by the received measurements,  our aim is to devise efficient algorithms for computing the WL centralized fusion estimators of the signal $\mathbf{x}(t)$, under the conditions of $\mathbb{T}_k$-properness, for $k=1,2$. 

With the purpose of a WL processing,  the $4n$-dimensional augmented vectors   are considered. Then, the centralized fusion estimation problem is addressed by applying the traditional estimation methods on the following WL  stacked state-space system:

\begin{eqnarray}
\bar{\mathbf{x}}(t+1)&=&  \bar{\boldsymbol{\Phi}}(t)\bar{\mathbf{x}}(t)+\bar{\mathbf{u}}(t), \ \ t\geq 0, \label{system3a} \\ 
\vec{\mathbf{z}}(t) &=&  \boldsymbol{\mathcal{C}} \bar{\mathbf{x}}(t)+\vec{\mathbf{v}}(t), \ \ t\geq 1, \label{system3b} \\ 
\vec{\mathbf{y}}(t) &=&  \bar{\boldsymbol{\mathcal{D}}}^{\vec{\boldsymbol{\gamma}}}(t) \vec{\mathbf{z}}(t) + \bar{\boldsymbol{\mathcal{D}}}^{(\boldsymbol{1}-\vec{\boldsymbol{\gamma}})}(t)  \vec{\mathbf{y}}(t-1), \ \ t\geq 2, \label{system3c}
\end{eqnarray}
with $\vec{\mathbf{y}}(1)=\vec{\mathbf{z}}(1)$, and where
 $\vec{\mathbf{z}}(t)=\left[ \bar{\mathbf{z}}^{(1)^{\texttt{T}}}(t),\ldots, \bar{\mathbf{z}}^{(R)^{\texttt{T}}}(t) \right]^{\texttt{T}}$, $\vec{\mathbf{v}}(t)=\left[ \bar{\mathbf{v}}^{(1)^{\texttt{T}}}(t),\ldots, \bar{\mathbf{v}}^{(R)^{\texttt{T}}}(t) \right]^{\texttt{T}}$, and $\vec{\mathbf{y}}(t)=\left[ \bar{\mathbf{y}}^{(1)^{\texttt{T}}}(t),\ldots, \bar{\mathbf{y}}^{(R)^{\texttt{T}}}(t) \right]^{\texttt{T}}$. 
Moreover,  \begin{equation*}
\bar{\boldsymbol{\Phi}}(t)=\left[
  \begin{array}{cccc}
\mathbf{F}_1(t) &\mathbf{F}_2(t) &\mathbf{F}_3(t) &\mathbf{F}_4(t) \\
\mathbf{F}_2^{*}(t) &\mathbf{F}_1^*(t) &\mathbf{F}_4^*(t) &\mathbf{F}_3^*(t) \\
    \mathbf{F}_3^{\eta}(t) &\mathbf{F}_4^{\eta}(t) & \mathbf{F}_1^{\eta}(t) & \mathbf{F}_2^{\eta}(t) \\
  \mathbf{F}_4^{\eta''}(t)&\mathbf{F}_3^{\eta''}(t) & \mathbf{F}_2^{\eta''}(t) &\mathbf{F}_1^{\eta''}(t) \\
  \end{array}
\right], 
\end{equation*}
$\bar{\boldsymbol{\mathcal{D}}}^{\vec{\boldsymbol{\gamma}}}(t) =\boldsymbol{\Upsilon}\diag(\vec{\boldsymbol{\gamma}}^{\mathrm{r}}(t))\boldsymbol{\Upsilon}^{\texttt{H}}$,  $ \bar{\boldsymbol{\mathcal{D}}}^{(\boldsymbol{1}-\vec{\boldsymbol{\gamma}})}(t) =\boldsymbol{\Upsilon}\diag\left(\boldsymbol{1}_{4nR}-\vec{\boldsymbol{\gamma}}^{\mathrm{r}}(t)\right)\boldsymbol{\Upsilon}^{\texttt{H}}$,   with~$\boldsymbol{\Upsilon}=\boldsymbol{I}_R\otimes \boldsymbol{\mathcal{T}}$ and ${\vec{\boldsymbol{\gamma}}}^r(t)=\left[{{\boldsymbol{\gamma}}^{(1)^{r^{\texttt{T}}}}}(t),\ldots,{{\boldsymbol{\gamma}}^{(R)^{r^{\texttt{T}}}}}(t)\right]^{\texttt{T}}$, and $\boldsymbol{\mathcal{C}}=\boldsymbol{1}_R\otimes \boldsymbol{I}_{4n}$. 

Furthermore, $\boldsymbol{\Gamma}_{\bar{\mathbf{u}}}(t,s)=\bar{\mathbf{Q}}(t) \delta_{t,s}$, $\boldsymbol{\Gamma}_{\vec{\mathbf{v}}}(t,s)=\vec{\mathbf{R}}(t)\delta_{t,s}$, and $\boldsymbol{\Gamma}_{\bar{\mathbf{u}}\vec{\mathbf{v}}}(t,s)=\vec{\mathbf{S}}(t)\delta_{t,s}$, where  $\vec{\mathbf{R}}(t)=\diag\left(\bar{\mathbf{R}}^{(1)}(t), \ldots, \bar{\mathbf{R}}^{(R)}(t)\right)$, with $\bar{\mathbf{R}}^{(i)}(t)=\boldsymbol{\Gamma}_{\bar{\mathbf{v}}^{(i)}}(t,t)$,  and $\vec{\mathbf{S}}(t)=\left[ \bar{\mathbf{S}}^{(1)}(t), \ldots, \bar{\mathbf{S}}^{(R)}(t) \right]$, with $\bar{\mathbf{S}}^{(i)}(t)=\boldsymbol{\Gamma}_{\bar{\mathbf{u}}\bar{\mathbf{v}}^{(i)}}(t,t)$, for $i=1,\dots, R$.

Now, the centralized fusion estimation problem is analyzed in a $\mathbb{T}_k$-properness setting. The following Proposition establishes conditions on system \eqref{system3a}-\eqref{system3c}  that guarantee the  $\mathbb{T}_k$-properness of the processes involved.  

\begin{Proposition}\label{prooo}
Given the WL stacked state-space model \eqref{system3a}-\eqref{system3c}, and taking into account the $\mathbb{T}_k$-properness concepts given in Definition 2,  the following properties can be established:
\begin{enumerate}
\item $\mathbf{x}(t)$ is $\mathbb{T}_1$-proper if and only if the initial state $\mathbf{x}(0)$ and the state noise $\mathbf{u}(t)$ are $\mathbb{T}_1$-proper, and the matrix $\bar{\boldsymbol{\Phi}}(t)$ is block diagonal as described below 
\begin{equation*}\label{F1}
\bar{\boldsymbol{\Phi}}(t)=\diag\left( \mathbf{F}_1(t), \mathbf{F}_1^{*}(t), \mathbf{F}_1^{\eta}(t), \mathbf{F}_1^{\eta''}(t)\right),
\end{equation*}


If additionally  $\mathbf{v}^{(i)}(t)$ is $\mathbb{T}_1$-proper, $\mathbf{u}(t)$  and $\mathbf{v}^{(i)}(t)$ are cross $\mathbb{T}_1$-proper,  and ${p}_{j,r}^{(i)}(t)={p}_{j,\eta}^{(i)}(t)={p}_{j,\eta'}^{(i)}(t)={p}_{j,\eta''}^{(i)}(t)\triangleq p_j^{(i)}(t)$, $\forall t$, $j=1,\dots,n$, $i=1,\dots, R$, then $\mathbf{x}(t)$ and $\mathbf{y}^{(i)}(t)$ are jointly $\mathbb{T}_1$-proper.

\item $\mathbf{x}(t)$ is $\mathbb{T}_2$-proper if and only if the initial state  $\mathbf{x}(0)$ and the state noise $\mathbf{u}(t)$ are $\mathbb{T}_2$-proper, and the matrix  $\bar{\boldsymbol{\Phi}}(t)$ is block diagonal as described below
\begin{equation*}\label{F2}
\bar{\boldsymbol{\Phi}}(t)=\diag\left( {\boldsymbol{\Phi}}_2(t), {\boldsymbol{\Phi}}_2^{\eta}(t)\right),\ \ \text{with}\ \ {\boldsymbol{\Phi}}_2(t)=\left[
  \begin{array}{cc}
 \mathbf{F}_1(t) & \mathbf{F}_2(t) \\   \mathbf{F}_2^{*}(t) & \mathbf{F}_1^{*}(t)
  \end{array}
\right],
\end{equation*}


If additionally $\mathbf{v}^{(i)}(t)$ is $\mathbb{T}_2$-proper, $\mathbf{u}(t)$  and $\mathbf{v}^{(i)}(t)$ are cross $\mathbb{T}_2$-proper, and ${p}_{j,r}^{(i)}(t)={p}_{j,\eta}^{(i)}(t)$, ${p}_{j,\eta'}^{(i)}(t)={p}_{j,\eta''}^{(i)}(t)$,  $\forall t$, $j=1,\dots,n$, $i=1,\dots, R$,  then $\mathbf{x}(t)$ and $\mathbf{y}^{(i)}(t)$ are jointly $\mathbb{T}_2$-proper.
\end{enumerate}
\end{Proposition}

\begin{Remark}
It should be observed that the conditions established in Proposition \ref{prooo} for ensuring the different type of properness on the processes involved in \eqref{system3a}-\eqref{system3c}, are similar to the one stated in \cite{Fernandez2021}.
\end{Remark}

Then,  under conditions of $\mathbb{T}_k$-properness, for $k=1,2$, the measurement equation  \eqref{system3c} in the above WL stacked state-space model can be expressed  in the following equivalent form of reduced dimension:
\begin{equation}\label{y_k}
{\mathbf{y}}_k(t)  =  \bar{\boldsymbol{\mathcal{D}}}_k^{\vec{\boldsymbol{\gamma}}}(t) \vec{\mathbf{z}}(t)+ 
     \bar{\boldsymbol{\mathcal{D}}}_k^{\mathbf{1}-\vec{\boldsymbol{\gamma}}} \vec{\mathbf{y}}(t-1),\quad t\geq 2,
\end{equation}
with ${\mathbf{y}}_k(1)=\boldsymbol{\Delta}_k\vec{\mathbf{z}}(1)$,  and $\boldsymbol{\Delta}_k=\mathbf{I}_R\otimes \left[\mathbf{I}_{kn},\mathbf{0}_{kn\times (4-k)n}\right]$.     Furthermore, $\bar{\boldsymbol{\mathcal{D}}}_k^{\vec{\boldsymbol{\gamma}}}(t)=\boldsymbol{\Upsilon}_k
\diag\left(\vec{\boldsymbol{\gamma}}^{r}(t)\right)\boldsymbol{\Upsilon}^{\texttt{H}}$ and  $\bar{\boldsymbol{\mathcal{D}}}_k^{\mathbf{1}-\vec{\boldsymbol{\gamma}}}(t)=\boldsymbol{\Upsilon}_k
\diag\left(\mathbf{1}_{4nR}-\vec{\boldsymbol{\gamma}}^{r}(t)\right)\boldsymbol{\Upsilon}^{\texttt{H}} $,  where~$\boldsymbol{\Upsilon}_k=\mathbf{I}_R\otimes\boldsymbol{\mathcal{T}}_k$, with  $\boldsymbol{\mathcal{T}}_k=\frac12\boldsymbol{\mathcal{B}}_k\otimes \mathbf{I}_n$,  and
   \begin{equation*}\label{Bk}
\boldsymbol{\mathcal{B}}_k=\left\{
\begin{array}{l}
      \left[1\quad \eta\quad \eta'\quad \eta''\right],\quad \text{for }k=1 \\
      \left[\begin{array}{cccc}1 & \eta & \eta'& \eta''\\
      1 & -\eta & \eta' & -\eta''\end{array}\right],\quad \text{for } k=2.    \end{array}\right.        \end{equation*}

In addition,
\begin{equation*}\label{Pi_k_vec}
 \begin{aligned}
 \bar{\boldsymbol{\Pi}}_k^{\vec{\boldsymbol{\gamma}}}(t) & =  E\left[\bar{\boldsymbol{\mathcal{D}}}_k^{\vec{\boldsymbol{\gamma}}}(t)\right]=\diag\left(\bar{\boldsymbol{\Pi}}_k^{\boldsymbol{\gamma}^{(1)}}(t),\ldots, \bar{\boldsymbol{\Pi}}_k^{\boldsymbol{\gamma}^{(R)}}(t)\right),\\
 \bar{\boldsymbol{\Pi}}_k^{(\boldsymbol{1}-\vec{\boldsymbol{\gamma}})}(t)& = E\left[\bar{\boldsymbol{\mathcal{D}}}_k^{(\boldsymbol{1}-\vec{\boldsymbol{\gamma}})}(t)\right]=\diag\left(\bar{\boldsymbol{\Pi}}_k^{\left(\boldsymbol{1}-\boldsymbol{\gamma}^{(1)}\right)}(t),\ldots, \bar{\boldsymbol{\Pi}}_k^{\left(\boldsymbol{1}-\boldsymbol{\gamma}^{(R)}\right)}(t)\right),
 \end{aligned}
        \end{equation*}
        with
$\bar{\boldsymbol{\Pi}}_k^{\boldsymbol{\gamma}^{(i)}}(t)=\left[\boldsymbol{\Pi}_k^{(i)}(t),\boldsymbol{0}_{kn\times(4-k)n}\right]$ and $\bar{\boldsymbol{\Pi}}_k^{\left(\boldsymbol{1}-\boldsymbol{\gamma}^{(i)}\right)}(t)=\left[\boldsymbol{I}_{kn}-\boldsymbol{\Pi}_k^{(i)}(t),\boldsymbol{0}_{kn\times(4-k)n}\right]$,  
\begin{equation}\label{Pi_k}
 \begin{aligned}
           \boldsymbol{\Pi}_1^{(i)}(t) & = \diag\left(p_{1,r}^{(i)}(t),\ldots,p_{n,r}^{(i)}(t)\right), \quad  i=1,\ldots,R,\\  
         \boldsymbol{\Pi}_2^{(i)}(t)  & =  \dfrac12\left[\begin{array}{cc}
                                                       \boldsymbol{\Pi}_{a}^{(i)}(t) & \boldsymbol{\Pi}_{b}^{(i)}(t) \\
                                                       \boldsymbol{\Pi}_{b}^{(i)}(t) & \boldsymbol{\Pi}_{a}^{(i)}(t)
                                                     \end{array}\right],\quad i=1,\ldots,R,
      \end{aligned}
        \end{equation}
        and
\begin{equation*}\label{Pi_2a,b}
      \begin{array}{c}
         \boldsymbol{\Pi}_{a}^{(i)}(t)=\diag\left(p_{1,r}^{(i)}(t)+p_{1,\eta'}^{(i)}(t),\ldots,p_{n,r}^{(i)}(t)+p_{n,\eta'}^{(i)}(t)\right),\quad i=1,\ldots,R,\\
         \boldsymbol{\Pi}_{b}^{(i)}(t)=\diag\left(p_{1,r}^{(i)}(t)-p_{1,\eta'}^{(i)}(t),\ldots,p_{n,r}^{(i)}(t)-p_{n,\eta'}^{(i)}(t)\right),\quad i=1,\ldots,R.        \end{array}
        \end{equation*}
        
\begin{Remark}\label{remark1}
It is worth to note  that $\mathbb{T}_k$-properness also allows to reduce the dimension of the equations \eqref{system3a} and \eqref{system3b} by replacing the $4n$-dimenstional augmented  processes $\bar{\mathbf{x}}(t)$, $\bar{\mathbf{u}}(t)$ $\bar{\mathbf{z}}^{(i)}(t)$, $\bar{\mathbf{v}}^{(i)}(t)$, and the matrix $\bar{\boldsymbol{\Phi}}(t)$ by the corresponding $kn$-dimensional vectors $\mathbf{x}_k(t)$, $\mathbf{u}_k(t)$, $\mathbf{z}^{(i)}_k(t)$, $\mathbf{v}^{(i)}_k(t)$ and $\boldsymbol{\Phi}_k(t)$, defined as
\begin{itemize}
 \item {\bf{$\mathbb{T}_1$-proper case:}} $\mathbf{x}_1(t)\triangleq \mathbf{x}(t)$, $\mathbf{u}_1(t)\triangleq \mathbf{u}(t)$, $\mathbf{z}^{(i)}_1(t)\triangleq \mathbf{z}^{(i)}(t)$, $\mathbf{v}^{(i)}_1(t)\triangleq \mathbf{v}^{(i)}(t)$ and $\boldsymbol{\Phi}_1(t)\triangleq \mathbf{F}_1(t)$.
  \item {\bf{$\mathbb{T}_2$-proper case:}} $\mathbf{x}_2(t)\triangleq \left[\mathbf{x}(t),\mathbf{x}^{\texttt{H}}(t)\right]^{\texttt{T}}$, $\mathbf{u}_2(t)\triangleq \left[\mathbf{u}(t),\mathbf{u}^{\texttt{H}}(t)\right]^{\texttt{T}}$, $\mathbf{z}^{(i)}_2(t)\triangleq \left[\mathbf{z}^{(i)}(t),\mathbf{z}^{(i)^{\texttt{H}}}(t)\right]^{\texttt{T}}$, $\mathbf{v}^{(i)}_2(t)\triangleq \left[\mathbf{v}^{(i)}(t),\mathbf{v}^{(i)^{\texttt{H}}}(t)\right]^{\texttt{T}}$, and $\boldsymbol{\Phi}_2(t)$ given in Proposition \ref{prooo},   in a $\mathbb{T}_2$-proper scenario.   
\end{itemize} 
 Furthermore,   
$\boldsymbol{\Gamma}_{\mathbf{u}_k}(t,s)=\mathbf{Q}_k (t) \delta_{t,s}$, $\boldsymbol{\Gamma}_{\mathbf{v}_k^{(i)}}(t,s)=\mathbf{R}_k^{(i)}(t) \delta_{t,s}$, $\boldsymbol{\Gamma}_{\mathbf{u}_k\mathbf{v}_k^{(i)}}(t,s)=\mathbf{S}_k^{(i)}(t) \delta_{t,s}$, and $\boldsymbol{\Gamma}_{\mathbf{x}_k}(0,0)=\mathbf{P}_{0_k}$.
\end{Remark}

Thus, whereas the optimal linear processing in the tessarine domain suggests computing the LLMS filter of the state  $\mathbf{x}(t)\in \mathbb{T}^n $ from its projection onto the augmented measurements $\{\vec{\mathbf{y}}(1),\ldots \vec{\mathbf{y}}(t)\}$,  under conditions of $\mathbb{T}_k$-properness, for $k=1,2$, this estimator can be obtained from the measurements $\{\mathbf{y}_k(1), \ldots, \mathbf{y}_k(t)\}$  defined in \eqref{y_k}, which gives rise to the so-called $\mathbb{T}_k$-proper estimators.  This approach supposes a reduction in the dimension of the problem that leads to computational savings that cannot be attained from a real formalism.  

This methodology has been recently applied to design recursive fusion estimation algorithms for multi-sensor systems affected by random delays and missing measurements \cite{Jimenez2021}.  In this paper, we are interested in extending this methodology to systems affected by random multiple-packet dropouts.

\section{\texorpdfstring{$\mathbb{T}_k$}-proper centralized fusion filtering estimation}

In this section, based on Kalman filter techniques, an efficient algorithm is provided for the computation of the $\mathbb{T}_k$-proper LLMS centralized fusion filter $ \hat{{\mathbf{x}}}^{\mathbb{T}_k}(t|t)$,  for $k=1,2$,  of the state $\mathbf{x}(t)$ described by the state-space system with packet dropouts given by the equations \eqref{system3a}-\eqref{system3b}, and \eqref{y_k}, 
as well as its associated error {\it pseudo} covariance matrix $ {\mathbf{P}}^{\mathbb{T}_k}(t|t)$. For this purpose, a recursive algorithm is devised under $ \mathbb{T}_k$-properness conditions for the projection of $\bar{\mathbf{x}}(t)$ onto the set of measurements $\{\mathbf{y}_k(1), \ldots, \mathbf{y}_k(t)\}$, denoted by $\hat{\mathbf{x}}_k(t|t)$, and its error {\it pseudo} covariance matrix $ {\mathbf{P}}_k(t|t)$. 
Then, $ \hat{{\mathbf{x}}}^{\mathbb{T}_k}(t|t)$ and ${\mathbf{P}}^{\mathbb{T}_k}(t|t)$ 
are determined by the first $n $ components of $\hat{{\mathbf{x}}}_k(t|t)$ and 
${\mathbf{P}}_k(t|t)$, respectively.

Theorem \ref{teo centralized filter} summarizes the formulas of this $\mathbb{T}_k$-proper LLMS centralized fusion filtering algorithm.


\begin{Theorem}
\label{teo centralized filter}
The $\mathbb{T}_k$-proper LLMS centralized fusion filter, $ \hat{{\mathbf{x}}}^{\mathbb{T}_k}(t|t)$, for $k=1,2$, is obtained as follows
\begin{equation*}
\begin{split}
\hat{{\mathbf{x}}}^{\mathbb{T}_1}(t|t)& =\hat{\mathbf{x}}_1(t|t),	\\ \hat{{\mathbf{x}}}^{\mathbb{T}_2}(t|t)& =\left[\boldsymbol{1}_n,\boldsymbol{0}_{n}\right]\hat{\mathbf{x}}_2(t|t),
\end{split}
\end{equation*}
where, for $k=1,2$, $\hat{{\mathbf{x}}}_k(t|t)$ is calculated from the recursive equation
\begin{equation}\label{filter}
   \hat{{\mathbf{x}}}_k(t|t)=\hat{{\mathbf{x}}}_k(t|t-1)+{\mathbf{L}}_k(t){\boldsymbol{\epsilon}}_k(t),\quad t\geq 1,
   \end{equation}
and $\hat{\mathbf{x}}_k(t+1|t)$ satisfies the recursive expression
\begin{equation}\label{pred}
  \hat{{\mathbf{x}}}_k(t+1|t)={\boldsymbol{\Phi}}_k(t)\hat{{\mathbf{x}}}_k(t|t)+{\mathbf{H}}_k(t){\boldsymbol{\epsilon}}_k(t), \quad t\geq 1,
\end{equation}
with initial values  $\hat{{\mathbf{x}}}_k(1|0)=\hat{{\mathbf{x}}}_k(0|0)=\mathbf{0}_{kn}$.    

The innovations ${\boldsymbol{\epsilon}}_k(t)$  are recursively calculated from the formula
\begin{equation}\label{central innovations}
  {\boldsymbol{\epsilon}}_k(t)=  {\mathbf{y}}_k(t)-\boldsymbol{\Pi}_k(t)\boldsymbol{\mathcal{C}}_k\hat{{\mathbf{x}}}_k(t|t-1)  -\left( \boldsymbol{I}_{knR}  - \boldsymbol{\Pi}_k(t)\right) \mathbf{y}_k(t-1),\quad t\geq 2,
\end{equation}
with initial value ${\boldsymbol{\epsilon}}_k(1)={\mathbf{y}}_k(1)$,  and $\boldsymbol{\mathcal{C}}_k=\mathbf{1}_R\otimes \mathbf{I}_{kn}$.

 Moreover,   ${\mathbf{H}}_k(t)={\mathbf{S}}_k(t){\boldsymbol{\Pi}}_k(t){\boldsymbol{\Omega}}_k^{{-1}}(t)$, where
${\mathbf{S}}_k(t)=[{\mathbf{S}}_k^{(1)}(t),\ldots,{\mathbf{S}}_k^{(R)}(t)]$, 
and  ${\boldsymbol{\Pi}}_k(t)=\diag\left({\boldsymbol{\Pi}}_k^{(1)}(t),\ldots,{\boldsymbol{\Pi}}_k^{(R)}(t)\right)$,   with ${\boldsymbol{\Pi}}_k^{(i)}(t)$ given in \eqref{Pi_k}.    ${\mathbf{L}}_k(t)={\boldsymbol{\Theta}}_k(t){\boldsymbol{\Omega}}_k^{{-1}}(t)$,      where 
 the matrices ${\boldsymbol{\Theta}}_k(t)$ is obtained from the equation
\begin{equation}\label{Theta}
\begin{aligned}
  {\boldsymbol{\Theta}}_k(t)&={\mathbf{P}}_k(t|t-1)\boldsymbol{\mathcal{C}}_k^{\texttt{T}}{\boldsymbol{\Pi}}_k(t),\quad t\geq 2;\\
     {\boldsymbol{\Theta}}_k(1)&= \mathbf{1}_R^{\texttt{T}}\otimes\mathbf{D}_{k}(1),
\end{aligned}
\end{equation}
with 
\begin{equation}\label{D_k(1)}
  \mathbf{D}_{k}(1)=\left[\mathbf{I}_{kn},\mathbf{0}_{kn\times(4-k)n}\right]{\boldsymbol{\Gamma}}_{\bar{\mathbf{x}}}(1,1)\left[\mathbf{I}_{kn},\mathbf{0}_{kn\times(4-k)n}\right]^{\texttt{T}},
\end{equation}
 and where ${\boldsymbol{\Gamma}}_{\bar{\mathbf{x}}}(t,t)$ is given by the recursive expression
 \begin{equation}\label{D local filter}
{\boldsymbol{\Gamma}}_{\bar{\mathbf{x}}}(t,t)=\bar{\boldsymbol{\Phi}}(t-1) {\boldsymbol{\Gamma}}_{\bar{\mathbf{x}}}(t-1,t-1)\bar{\boldsymbol{\Phi}}^{\texttt{H}}(t-1)+\bar{\mathbf{Q}}(t-1),\quad t\geq 1;\quad {\boldsymbol{\Gamma}}_{\bar{\mathbf{x}}}(0,0)=\bar{\mathbf{P}}_0.  
\end{equation}
In addition, 
\begin{equation}\label{Omega central filter theorem}
\begin{aligned}
{\boldsymbol{\Omega}}_k(t)&= \boldsymbol{\Upsilon}_k\left\{\Cov\left(\vec{\boldsymbol{\gamma}}^r(t)\right)\circ \left({\boldsymbol{\Psi}}_{1}(t)-{\boldsymbol{\Psi}}_{2}(t)- {\boldsymbol{\Psi}}_{2}^{{\texttt{H}}}(t)+{\boldsymbol{\Psi}}_{3}(t)\right)\right\} \boldsymbol{\Upsilon}_k^{\texttt{H}}
\\ & \quad + \boldsymbol{\Upsilon}_k\left\{\boldsymbol{\Gamma}_{\vec{\boldsymbol{\gamma}}^r}(t,t)\circ
\left(\boldsymbol{\Upsilon}^{\texttt{H}}\vec{\mathbf{R}}(t)\boldsymbol{\Upsilon}\right)\right\} \boldsymbol{\Upsilon}_k^{\texttt{H}}
+  {\boldsymbol{\Pi}}_k(t)\boldsymbol{\mathcal{C}}_k{\mathbf{P}}_k(t|t-1)\boldsymbol{\mathcal{C}}_k^{\texttt{T}}{\boldsymbol{\Pi}}_k(t),
\end{aligned}
\end{equation}
where
 \begin{equation*}
\begin{aligned}
       {\boldsymbol{\Psi}}_{1}(t)& =      
\boldsymbol{\Upsilon}^{\texttt{H}}\boldsymbol{\mathcal{C}}{\boldsymbol{\Gamma}}_{\bar{\mathbf{x}}}(t,t)\boldsymbol{\mathcal{C}}^{\texttt{T}}\boldsymbol{\Upsilon},\\
{\boldsymbol{\Psi}}_{2}(t)& = 
\boldsymbol{\Upsilon}^{\texttt{H}}\boldsymbol{\mathcal{C}}\left( \bar{\boldsymbol{\Phi}}(t-1){\boldsymbol{\Gamma}}_{\bar{\mathbf{x}}\vec{\mathbf{y}}}(t-1,t-1) + \vec{\mathbf{S}}(t-1) \bar{\boldsymbol{\Pi}}^{\vec{\gamma}}(t-1)  \right) \boldsymbol{\Upsilon},\\
{\boldsymbol{\Psi}}_{3}(t)&=
\boldsymbol{\Upsilon}^{\texttt{H}}{\boldsymbol{\Gamma}}_{\vec{\mathbf{y}}}(t-1, t-1)\boldsymbol{\Upsilon},\\
     \end{aligned}\end{equation*}
with $\boldsymbol{\Gamma}_{\bar{\mathbf{x}}}(t,t)$ computed in (\ref{D local filter}), 
\begin{equation*}\label{Dxy}
\begin{aligned}
{\boldsymbol{\Gamma}}_{\bar{\mathbf{x}}\vec{\mathbf{y}}}(t,t)&={\boldsymbol{\Gamma}}_{\bar{\mathbf{x}}}(t,t)\boldsymbol{\mathcal{C}}^{\texttt{T}} \bar{\boldsymbol{\Pi}}^{\vec{\gamma}}(t)+\left( \bar{\boldsymbol{\Phi}}(t-1){\boldsymbol{\Gamma}}_{\bar{\mathbf{x}}\vec{\mathbf{y}}}(t-1,t-1) + \vec{\mathbf{S}}(t-1) \bar{\boldsymbol{\Pi}}^{\vec{\gamma}}(t-1)  \right)\bar{\boldsymbol{\Pi}}^{\boldsymbol{1}-\vec{\gamma}}(t),\quad t \geq 2;\\
{\boldsymbol{\Gamma}}_{\bar{\mathbf{x}}\vec{\mathbf{y}}}(1,1)&= {\boldsymbol{\Gamma}}_{\bar{\mathbf{x}}}(t,t)\boldsymbol{\mathcal{C}}^{\texttt{T}}, 
\end{aligned}
     \end{equation*}
     and 
     \begin{equation*}\label{Dy}
     \begin{split}
{\boldsymbol{\Gamma}}_{\vec{\mathbf{y}}}(t,t)&=\boldsymbol{\Upsilon}\left\{\boldsymbol{\Gamma}_{\vec{\boldsymbol{\gamma}}^r}(t,t) \circ \left( {\boldsymbol{\Psi}}_{1}(t)+\boldsymbol{\Upsilon}^{\texttt{H}}\vec{\mathbf{R}}(t)\boldsymbol{\Upsilon} \right)+  \boldsymbol{\Gamma}_{\vec{\boldsymbol{\gamma}}^r (\boldsymbol{1}-\vec{\boldsymbol{\gamma}}^r)}(t,t) \circ {\boldsymbol{\Psi}}_{2}(t)  +  \boldsymbol{\Gamma}_{\vec{\boldsymbol{\gamma}}^r (\boldsymbol{1}-\vec{\boldsymbol{\gamma}}^r)}^{\texttt{T}}(t,t) \circ {\boldsymbol{\Psi}}_{2}^{\texttt{H}}(t)  \right.  \\  & \qquad \qquad \left.  +
 \boldsymbol{\Gamma}_{(\boldsymbol{1}-\vec{\boldsymbol{\gamma}}^r)}(t,t) \circ {\boldsymbol{\Psi}}_{3}(t) \right\}\boldsymbol{\Upsilon}^{\texttt{H}},\quad t \geq 2;\\
{\boldsymbol{\Gamma}}_{\vec{\mathbf{y}}}(1,1)&= \boldsymbol{\mathcal{C}}{\boldsymbol{\Gamma}}_{\bar{\mathbf{x}}}(t,t)\boldsymbol{\mathcal{C}}^{\texttt{T}}+\vec{\mathbf{R}}(1).
 \end{split}
    \end{equation*}

Finally,  the $\mathbb{T}_k$-proper centralized fusion filtering error {\it pseudo} covariance matrix, ${\mathbf{P}}^{\mathbb{T}_k}(t|t)$, for $k=1,2$, is obtained 
as follows:
\begin{equation*}
\begin{split}
{\mathbf{P}}^{\mathbb{T}_1}(t|t)&={\mathbf{P}}_1(t|t),\\
{\mathbf{P}}^{\mathbb{T}_2}(t|t)&=\left[\boldsymbol{1}_n,\boldsymbol{0}_n\right]{\mathbf{P}}_2(t|t)\left[\boldsymbol{1}_n,\boldsymbol{0}_n\right]^{\texttt{T}},
\end{split}
\end{equation*}
where, for $k=1,2$, ${\mathbf{P}}_k(t|t)$ satisifes 
the following recursive equation:
\begin{equation}\label{centralerrorfil}
  {\mathbf{P}}_k(t|t)={\mathbf{P}}_k(t|t-1)-{\boldsymbol{\Theta}}_k(t){\boldsymbol{\Omega}}_k^{{-1}}(t){\boldsymbol{\Theta}}_k^{{\texttt{H}}}(t),
\end{equation}
with  initial condition ${\mathbf{P}}_k(0|0)={\mathbf{P}}_{0_k}$, and 
\begin{equation}\label{centralerrorpred}
\begin{aligned}
{\mathbf{P}}_k(t+1|t)&={\boldsymbol{\Phi}}_k(t){\mathbf{P}}_k(t|t){\boldsymbol{\Phi}}_k^{\texttt{H}}(t)
  -{\boldsymbol{\Phi}}_k(t){\boldsymbol{\Theta}}_k(t){\mathbf{H}}_k^{{\texttt{H}}}(t)-{\mathbf{H}}_k(t){\boldsymbol{\Theta}}_k^{{\texttt{H}}}(t){\boldsymbol{\Phi}}_k^{\texttt{H}}(t)\\
  &\qquad
  -{\mathbf{H}}_k(t){\boldsymbol{\Omega}}_k(t){\mathbf{H}}_k^{{\texttt{H}}}(t)+{\mathbf{Q}}_k(t),
\end{aligned}
\end{equation}
with  initial condition ${\mathbf{P}}_k(1|0)=\mathbf{D}_{k}(1)$ .    
\end{Theorem}

\begin{Remark}\label{r6}
{Notice that the computational load of the $\mathbb{T}_k$-proper LLMS centralized fusion filtering algorithms, for $k=1,2$, given in Theorem \ref{teo centralized filter} is the same as that of their quaternion domain  counterparts,  i.e., those derived by using a quaternion strictly linear (QSL) and quaternion semi-widely linear (QSWL) processing, respectively. 

As a consequence, it is noteworthy to see that the proposed $\mathbb{T}_k$-proper LLMS centralized fusion filtering algorithm provides estimations of the state that is equivalent to the one obtained from a WL processing or a real vectorial processing, whereas the computational load implied is reduced from $\mathcal{O}(64R^3n^3)$ to $\mathcal{O}(kR^3n^3)$, for $k=1,2$ \cite{Nitta2019}.}

\end{Remark}

\section{Numerical Example}
Our aim in this  section is to numerically analyze the performance and benefits of the $\mathbb{T}_k$-proper LLMS centralized fusion filtering algorithm proposed in Theorem \ref{teo centralized filter}. Two examples are proposed: the first one, from simulated values, in which  a scalar signal is estimated from the observations provided by several sensors; and the second one, a realistic model of a bidimensional tessarine state-space model which described a great amount of experimental phenomena. In both examples, by varying the Bernoulli parameters, different situations are compared in order to illustrate the effectiveness of the proposed algorithm in both $\mathbb{T}_k$-proper scenarios, for $k=1,2$. 
 

\subsection{Example 1}

Consider the following  multi-sensor tessarine state-space system: 
\begin{equation}\label{available obs equation in example}
\begin{split}
x(t+1)&=F_1(t)x(t)+u(t),\quad t\geq 0,\\
z^{(i)}(t)&=x(t)+v^{(i)}(t),\quad t\geq 1,\\
y^{(i)}(t)&=\gamma^{(i)}(t)\star z^{(i)}(t)+(1-\gamma^{(i)}(t))\star y^{(i)}(t-1),\quad t\geq 2,
\end{split}
\end{equation} 
%
%
for $i=1,\dots,R$, with $ y^{(i)}(1)=z^{(i)}(1)$, and where  $F_1(t)=0.3+0.3\eta+0.1\eta'+0.2\eta''\in\mathbb{T}$.   Moreover,  $u(t)$ is a tessarine noise such that the covariance matrix of the associated real vector $\mathbf{u}^{r}(t)$ is of the form
\begin{equation} \boldsymbol{\Gamma}_{\mathbf{u}^{r}}(t,s)=\left(\begin{array}{cccc}
                        a & 0 & c & 0 \\
                        0 & b & 0 & c \\
                        c & 0 & a & 0 \\
                        0 & c & 0 & b
                      \end{array}\right)\delta_{t,s},
                      \label{Q in the example}
                      \end{equation}       
where the parameters $a$, $b$ and $c$ take different values depending on the $\mathbb{T}_k$-proper scenario considered. Furthermore, to guarantee the correlation hypothesis between the state and observation noises, $u(t)$ and $v^{(i)}(t)$,  they are assumed to satisfies the following expression: 
$$v^{(i)}(t)=\alpha_i u(t)+w^{(i)}(t),\quad t\geq 1,$$
with $\alpha_i\in \mathbb{R}$, and where, at each $i$, $w^{(i)}(t)$ is a tessarine white Gaussian noise independent of $u(t)$, whose real covariance matrix is given by
$$\boldsymbol{\Gamma}_{\mathbf{w}^{(i)^r}}(t,s)=\diag\left(\beta_i,\beta_i,\beta_i,\beta_i\right), \quad t\geq 1,$$
Specifically,  the following values of $\alpha_i$ and $\beta_i$, for $i=1,2,3,4,5$, will be considered in our simulations:
$$\alpha_1=0.5,\ \alpha_2=0.3,\ \alpha_3=0.9,\ \alpha_4=0.6,\ \alpha_5=0.2$$
$$\beta_1=95,\ \beta_2=125,\ \beta_3=87,\ \beta_4=83,\ \beta_5=73$$

Additionally,  the variance matrix of the real initial state $\mathbf{x}^r(0)$ is assumed to be of the form
\begin{equation} 
  \boldsymbol{\Gamma}_{\mathbf{x}^r}(0,0)=\left(\begin{array}{cccc}
                        d & 0 & f & 0 \\
                        0 & e & 0 & f \\
                        f & 0 & d & 0 \\
                        0 & f & 0 & e
                      \end{array}\right),
\label{P0 in the example}\end{equation}                      
whose values $d$, $e$ and $f$ will be specified in Section \ref{subT1} and Section \ref{subT2}, according to the different $\mathbb{T}_k$-proper scenario analyzed.

\subsubsection{\texorpdfstring{$\mathbb{T}_1$}-proper scenario}\label{subT1}
To guarantee that $x(t)$ and $y^{(i)}(t)$ are joint $\mathbb{T}_1$-proper, it has been taken $a=b=1$, $c=-0.5$ in \eqref{Q in the example} and $d=e=4$, $c=1.5$ in \eqref{P0 in the example}. 
Moreover, it has also been assumed that the components of the multiplicative noise in \eqref{available obs equation in example}, $\gamma_{\nu}^{(i)}(t)$ have constant probabilities  
$p_{\nu}^{(i)}=p^{(i)}$, for all $\nu=r,\eta,\eta'\eta''$, and $i=1,\ldots,R$.



\begin{figure}[ht]
\centering
\includegraphics[width=14cm]{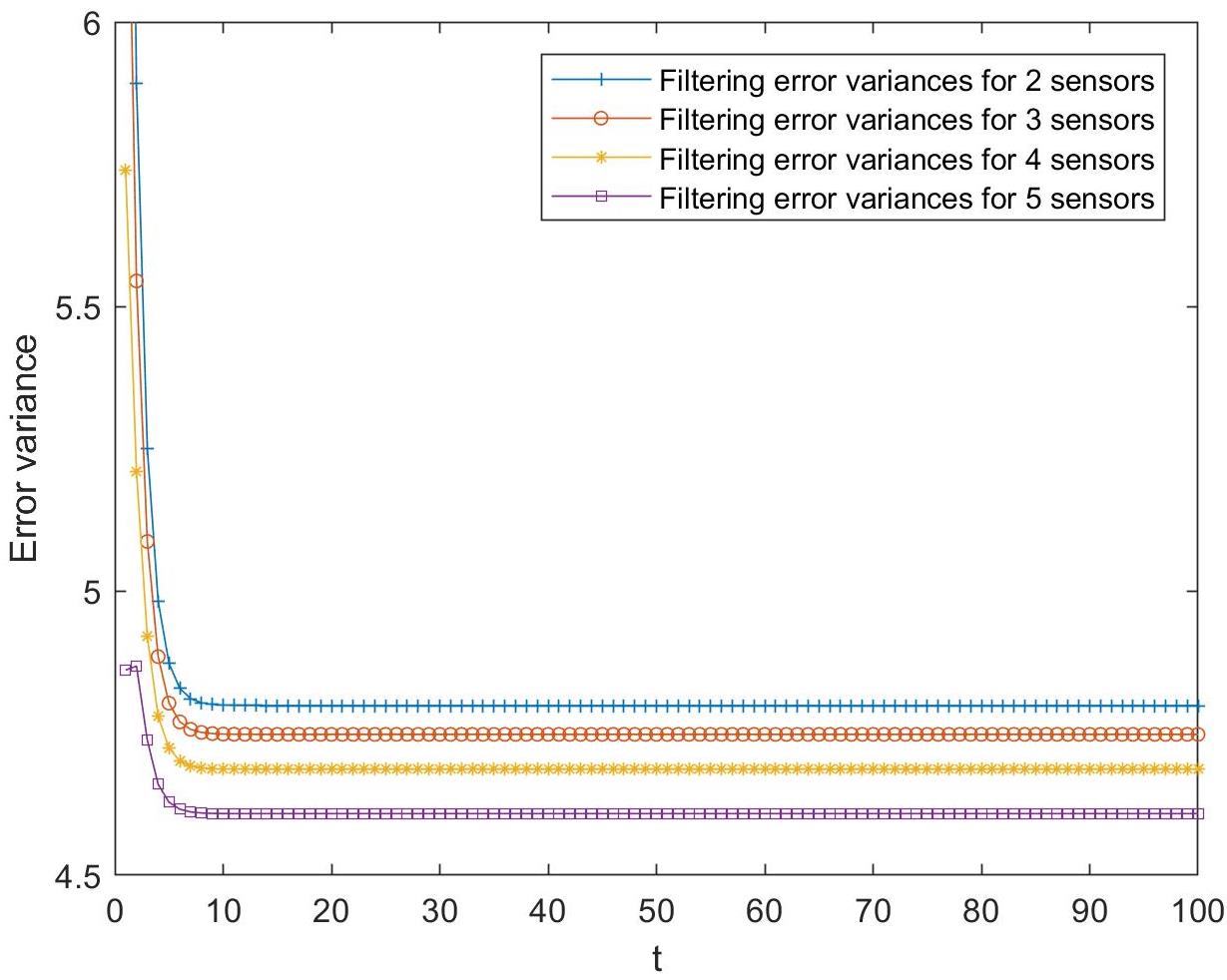}
\caption{\texorpdfstring{$\mathbb{T}_1$}-proper centralized fusion filtering error variances with 2, 3, 4, and 5 sensors.\label{fig2}}
\end{figure}

Firstly, the behavior of the estimators proposed is analyzed by considering a different number of sensors. Specifically, Figure \ref{fig2} shows the $\mathbb{T}_1$-proper centralized fusion filtering error variances computed from the observations provided by 2, 3, 4, and 5 sensors.  As expected, it is the estimators perform better  as the number of sensors increases, which makes sense because the number of observations used to estimate the signal increases.     

Next, in order to show the computational savings attained with the solution proposed, under $\mathbb{T}_1$-properness conditions, Table \ref{tab1} presents the computation time required to apply the  $\mathbb{T}_1$-proper centralized fusion filtering algorithms given in Theorem \ref{teo centralized filter}, and the conventional one devised from a real-valued linear processing  in the cases of $2$, $3$, $4$, and $5$ sensors. Then, a reduction in the computation time can be observed when the methodology proposed is used,  and this computational saving becomes more significant as the number of sensors increases.     

\begin{table}[ht]
\caption{Computation time (in seconds) for the \texorpdfstring{$\mathbb{T}_1$}- proper processing and the conventional one.\label{tab1}}
	\begin{tabular}{ccccccc}
			\toprule
	\multicolumn{1}{c}{ \textbf{Type of processing}}	& 		 & \multicolumn{4}{c}{ \textbf{Number of sensors}}  & 	     \\
			\cmidrule{3-6}
			 & &   $2$ & $3$ & $4$ & $5$  \\
			\midrule 
		\hskip-0.2cm	\begin{tabular}{c} $\mathbb{T}_1$- proper\\ Real-valued \end{tabular} & & \begin{tabular}{c} 4.552597\\ 5.367786 \end{tabular} & \begin{tabular}{c} 9.328620\\ 10.935009\end{tabular} & \begin{tabular}{c} 16.299967 \\ 18.187468\end{tabular} & \begin{tabular}{c} 25.112891\\ 27.570617\end{tabular} 
			\\ 
			\bottomrule
						 			\end{tabular}
\end{table}

\begin{figure}[ht]
\centering
\includegraphics[width=14cm]{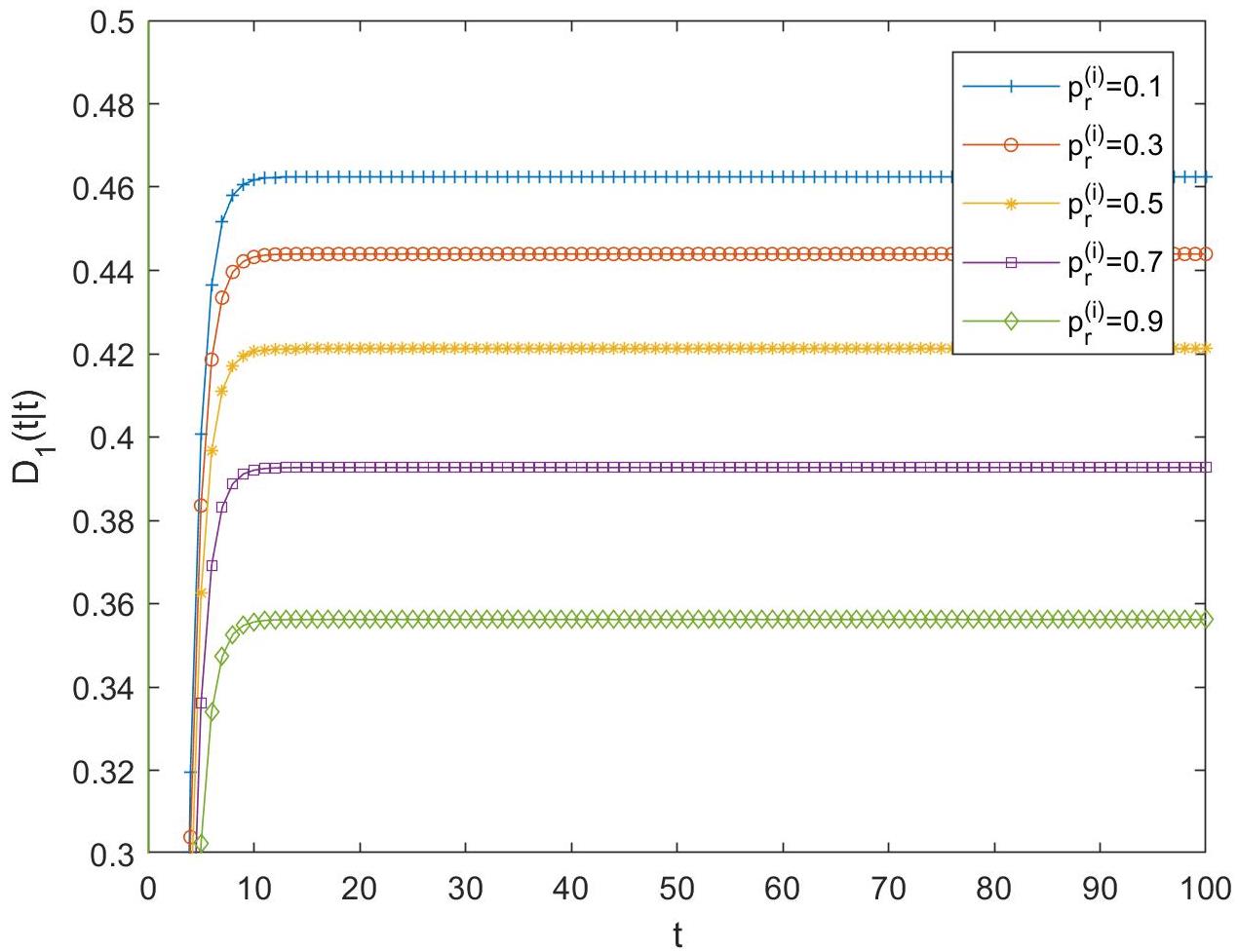}
\caption{Difference $D_1(t|t)$ between QSL and \texorpdfstring{$\mathbb{T}_1$}-proper centralized fusion filtering error variances for cases 1, 2, 3, 4,  and 5.\label{fig3}}
\end{figure}

Our second objective is to compare  tessarine and quaternion signal processing for different probabilities of updated/missing observations, under $\mathbb{T}_1$-properness conditions. For this purpose, the error variances of both $\mathbb{T}_1$ and QSL centralized fusion filters  have been calculated for the following cases:   
\begin{itemize}
      \item [-] Case 1: $p^{(i)}=0.1$, $\forall i=1,\ldots,5$;
      \item [-] Case 2: $p^{(i)}=0.3$, $\forall i=1,\ldots,5$;
      \item [-] Case 3: $p^{(i)}=0.5$, $\forall i=1,\ldots,5$;
      \item [-] Case 4: $p^{(i)}=0.7$, $\forall i=1,\ldots,5$;
      \item [-] Case 5: $p^{(i)}=0.9$, $\forall i=1,\ldots,5$.
\end{itemize}
Then, the difference between both tessarine and quaternion LLMS centralized fusion filtering error  variances, that is, $D_1(t|t)=P_{QSL}(t|t)-P_1(t|t)$, have been computed and displayed in Figure \ref{fig3}. In this figure, positive differences can be observed in all the cases, meaning that  it can be noted that  $\mathbb{T}_1$-proper  fusion estimators perform better than their quaternion counterparts. As expected, the fact that the $\mathbb{T}_1$-properness conditions are satisfied, determine that it is more appropriate to use
the $\mathbb{T}_1$-proper signal processing than the quaternion one, since it yields to better estimations. Moreover, these differences become smaller as the probability of updated observations increases.         

\begin{figure}[ht]
\centering
\includegraphics[width=14cm]{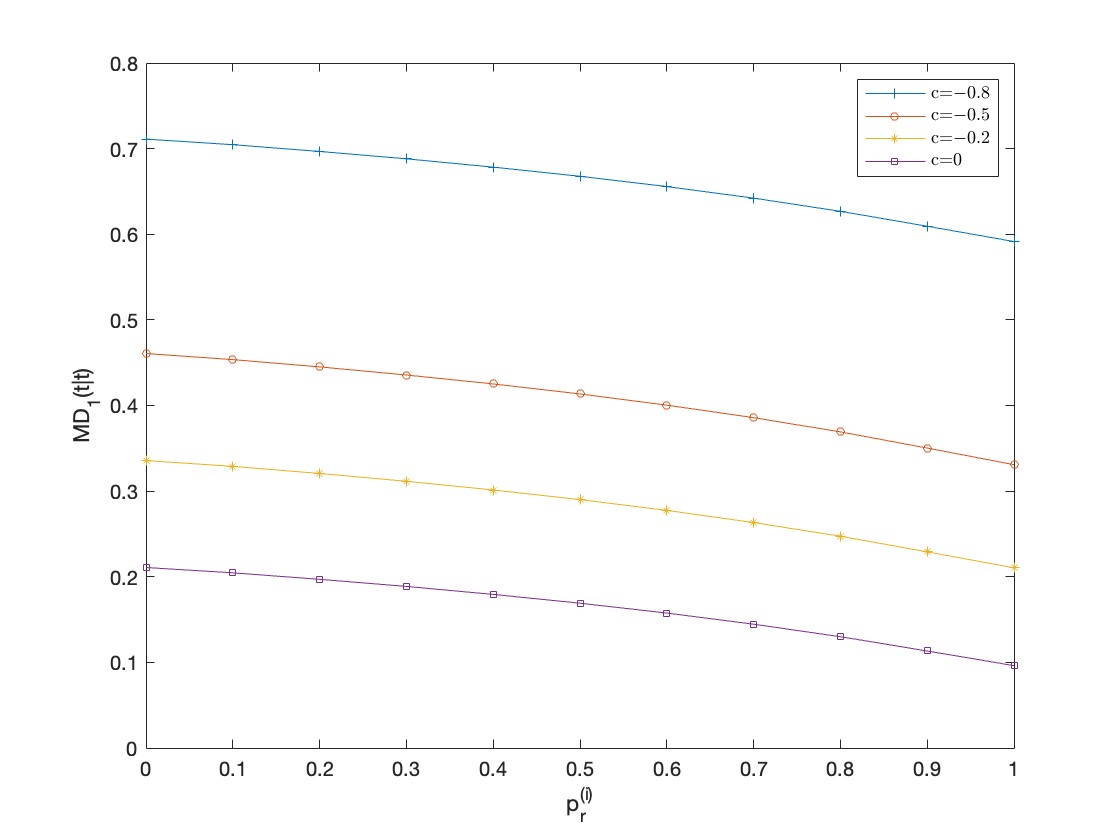}
\caption{Mean of the differences $MD_1(t|t)$ between QSL and \texorpdfstring{$\mathbb{T}_1$}-proper centralized fusion filtering error variances.\label{fig4}}
\end{figure}

Finally, with the aim of comparing both $QSL$ and $\mathbb{T}_1$-proper signal processing, they are applied by taking a fixed value for the probabilities of the Bernoulli parameters in all the sensors, but different values of $c$ in \eqref{Q in the example}, that is, $c=-0.8, -0.5, -0.2, 0$. Note that, for $c=0$, the state additive noise, $u(t)$, is $\mathbb{T}_1$ besides $\mathbb{Q}$-proper, and as $c$ is further away from 0, the $\mathbb{Q}$-properness conditions are further away. In this setting, the error variances of both $\mathbb{T}_1$-proper and QSL LLMS centralized fusion filters have been computed, and the mean of the differences between them, $MD_1(t|t)=\text{mean}(D_1(t|t))$,  have been displayed in Figure \ref{fig4} for the different values of $c$. In this figure, tessarine estimators are shown to be more accurate the further the noise $u(t)$ is from the $\mathbb{Q}$-properness conditions. Moreover, as in Figure \ref{fig3}, these differences decrease as the probability of updated  observations increases.

\subsubsection{\texorpdfstring{$\mathbb{T}_2$}-proper scenario}\label{subT2}

In order to guarantee $\mathbb{T}_2$-properness conditions,  $a=1$, $b=2$, $c=-0.5$ in \eqref{Q in the example}  is assumed, and $d=4$, $e=3$, $c=1.5$ in \eqref{P0 in the example}, and also $p_{r}^{(i)}=p_{\eta}^{(i)}$, and $p_{\eta'}^{(i)}=p_{\eta''}^{(i)}$,  for  $i=1,\ldots,5$.
 


\begin{figure}[ht]
\centering
\includegraphics[width=14cm]{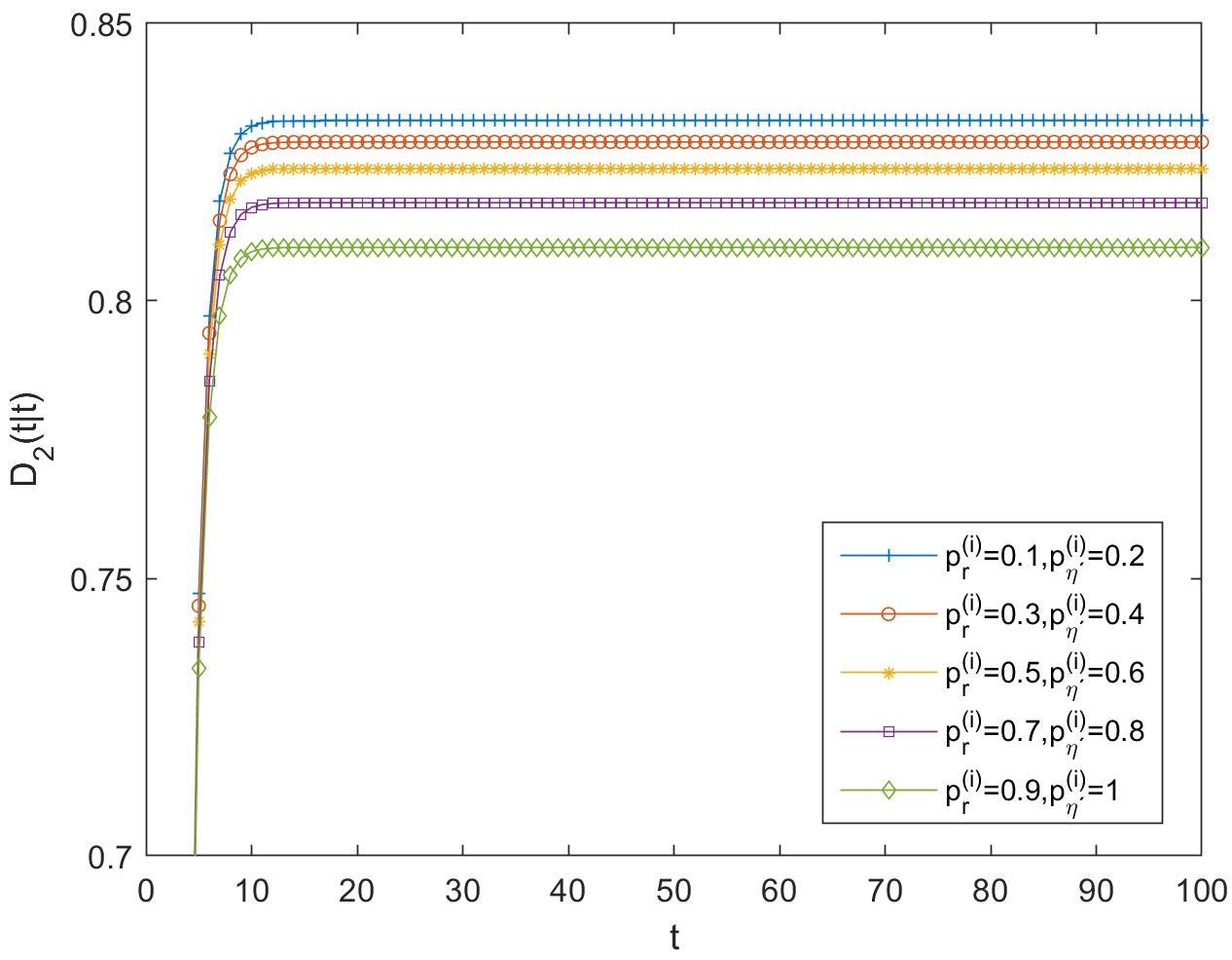}
\caption{Difference $D_2(t|t)$ between QSWL and \texorpdfstring{$\mathbb{T}_2$}-proper centralized fusion filtering error variances for cases 6, 7, 8, 9, and 10.\label{fig6}}
\end{figure}

As in the previous subsection, in order to compare the performance of tessarine and quaternion processing under $\mathbb{T}_2$-properness conditions, the differences between the LLMS centralized fusion filtering error variances of the QSWL and $\mathbb{T}_2$-proper estimators, denoted by $D_2(t|t)$, have been computed and displayed in Figure \ref{fig6}, for the following cases:
 \begin{itemize}
      \item [-] Case 6: $p_{r}^{(i)}=0.1$ and $p_{\eta'}^{(i)}=0.2$, 
      $\forall i=1,\ldots,5$;
      \item [-] Case 7: $p_{r}^{(i)}=0.3$ and $p_{\eta'}^{(i)}=0.4 $, 
      $\forall i=1,\ldots,5$;
      \item [-] Case 8: $p_{r}^{(i)}=0.5$ and $p_{\eta'}^{(i)}=0.6 $, 
      $\forall i=1,\ldots,5$;
      \item [-] Case 9: $p_{r}^{(i)}=0.7$ and $p_{\eta'}^{(i)}=0.8 $, 
      $\forall i=1,\ldots,5$;
      \item [-] Case 10: $p_{r}^{(i)}=0.9$ and $p_{\eta'}^{(i)}=1 $, 
      $\forall i=1,\ldots,5$.
      \end{itemize}
 Because these differences are positive in all the cases, the superiority of the $\mathbb{T}_2$-proper tessarine processing over the QSWL processing  under $\mathbb{T}_2$-properness conditions is clear, and these differences become smaller as the probability that the components of the available observation are updated increases.    


\subsection{Example 2}

Let us consider the following general equation of motion \cite{Talebi2020}:
\begin{equation}\label{motionEq}
\frac{\partial \varphi}{\partial t}=\phi, \qquad \     \text{and} \qquad   \frac{\partial \phi}{\partial t}=\upsilon, 
\end{equation}
where $ \varphi$ is the variable of interest, $\phi $ its range of change, and $\upsilon$ the input of the system.

Notice that equation (\ref{motionEq}) models a great amount of physical phenomena and it has been used, for example, in bearing-only tracking applications and rotation tracking problems, where $\upsilon$   represents,  respectively, the force or acceleration and the torque or angular acceleration.


In discrete-time, by taking $\mathbf{x}(t)=\left[ \varphi(t), \phi (t)\right]^{\texttt{T}}$, it is possible to build a model equivalent to that given in (\ref{motionEq}), as follows:
\begin{equation*}\label{motionEq2}
\begin{aligned}
\mathbf{x}(t+1)&=\left( \begin{array}{cc}   1 & 0.04 \\ 0& 1   \end{array}\right)  \mathbf{x}(t) + \left[\begin{array}{c} 0.0008 \\ 0.04 \end{array} \right] \varpi (t), \quad   t=1,\dots, 100;\\
\mathbf{x}(0)&=\boldsymbol{0}_{2\times 1},
  \end{aligned}
\end{equation*}
 where $\varpi (t)$ is a tessarine white noise with real covariance matrix:
 \begin{equation*}
E\left[  \varpi^r(t) \varpi^{r^{\texttt{T}}} (s)  \right]=\left( \begin{array}{cccc}  3& 0& 2 & 0\\ 0& 3 & 0 & 2\\  2 & 0 & 3 & 0\\  0 & 2  & 0& 3    \end{array}  \right) \delta_{ts},\quad t,s=1,\dots, 100.
\end{equation*}
  
  Moreover, the additive noise of single-sensor real observation equation, $\mathbf{v}(t)=\left[v_1(t), v_2(t)\right]^{\texttt{T}}$,  is assumed to be a tessarine white noise  with independent components and associated real covariance matrices given by 
\begin{equation*}
E\left[ \mathbf{v}_j^r(t) \mathbf{v}_j^{r^{\texttt{T}}} (s)  \right]=\left( \begin{array}{cccc}  6.5& 0& 0.1 & 0\\ 0& 6.5 & 0 &0.1\\  0.1 & 0 & 6.5 & 0\\  0 & 0.1  & 0& 6.5    \end{array}  \right) \delta_{ts},\quad t,s=1,\dots, 100 ,\  j=1,2.
\end{equation*} 

In order to guarantee the $\mathbb{T}_k$-properness conditions, the following assumptions and cases about the parameters of the Bernoulli random variables have been considered:
\begin{itemize}
      \item [$\bullet$] In the $\mathbb{T}_1$-proper scenario: $p_{j,\nu}(t)=p_{j}$, for all $j=1,2$, $\nu=r,\eta,\eta'\eta''$:
\begin{itemize}
      \item [-] Case 11: $p_1=0.1$, $p_2=0.2$;
      \item [-] Case 12: $p_1=0.3$, $p_2=0.4$;
      \item [-] Case 13: $p_1=0.5$, $p_2=0.6$;
      \item [-] Case 14: $p_1=0.7$, $p_2=0.8$;
      \item [-] Case 15: $p_1=0.9$, $p_2=1$.
\end{itemize}
\item [$\bullet$] In the $\mathbb{T}_2$-proper scenario: $p_{j,r}(t)=p_{j,\nu}(t)=p_{j,r}$, $p_{j,\eta'}(t)=p_{j,\nu''}(t)=p_{j,\eta'}$, for all $j=1,2$:
\begin{itemize}
      \item [-] Case 16: $p_{1,r}=0.1$, $p_{1,\eta'}=p_{2,r}=0.2$, $p_{2,\eta'}=0.3$;
      \item [-] Case 17: $p_{1,r}=0.3$, $p_{1,\eta'}=p_{2,r}=0.4$, $p_{2,\eta'}=0.5$;
      \item [-] Case 18: $p_{1,r}=0.5$, $p_{1,\eta'}=p_{2,r}=0.6$, $p_{2,\eta'}=0.7$;
      \item [-] Case 19: $p_{1,r}=0.7$, $p_{1,\eta'}=p_{2,r}=0.8$, $p_{2,\eta'}=0.9$;
      \item [-] Case 20: $p_{1,r}=0.9$, $p_{1,\eta'}=p_{2,r}=0.95$, $p_{2,\eta'}=1$.
\end{itemize}
\end{itemize}

\begin{figure}[ht]
\centering
\includegraphics[width=14cm]{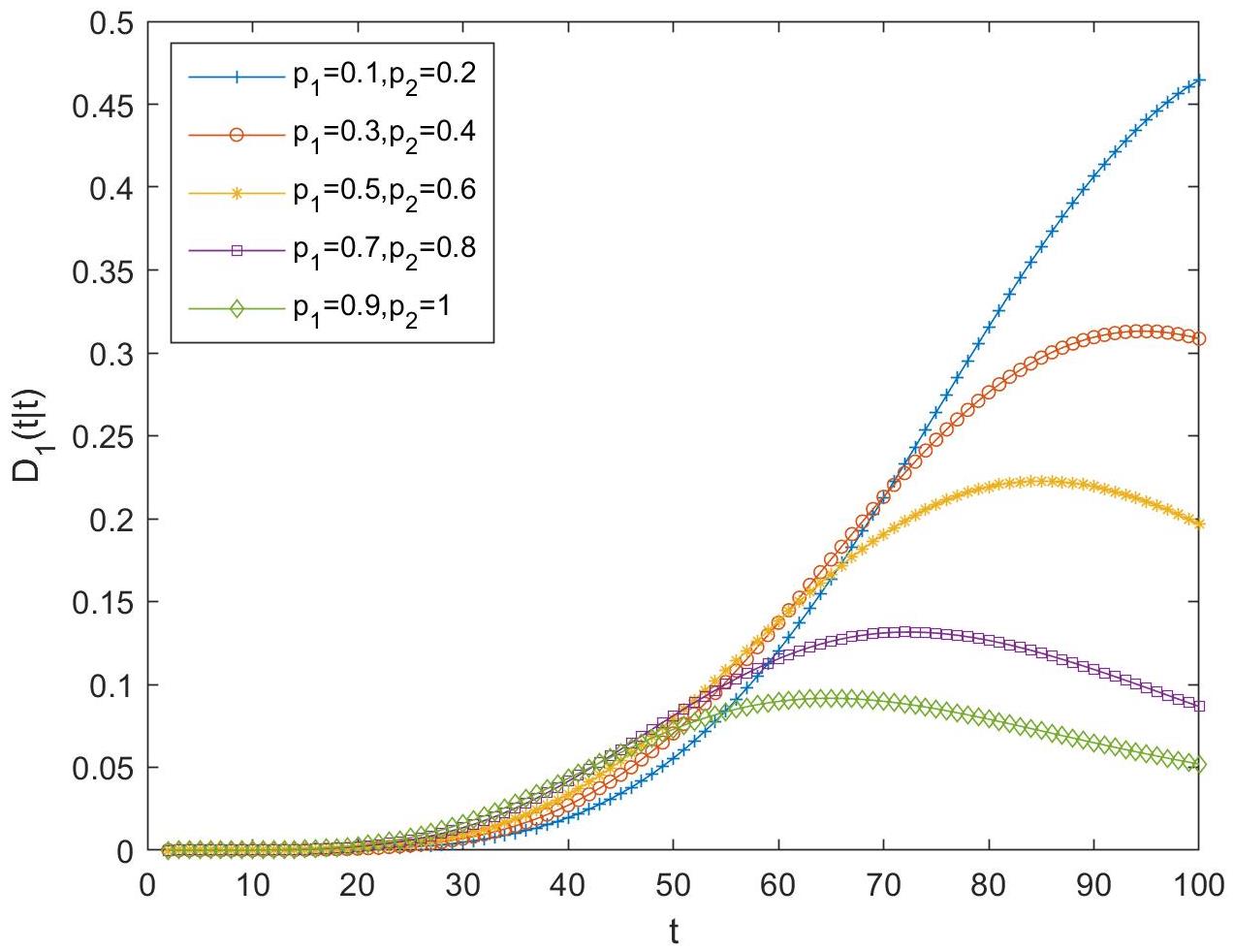}
\caption{Difference $D_1(t|t)$ between QSL and \texorpdfstring{$\mathbb{T}_1$}-proper fusion filtering error variances for the first component of the signal for cases 11, 12, 13, 14, and 15.\label{fig7}}
\end{figure}

\begin{figure}[ht]
\centering
\includegraphics[width=14cm]{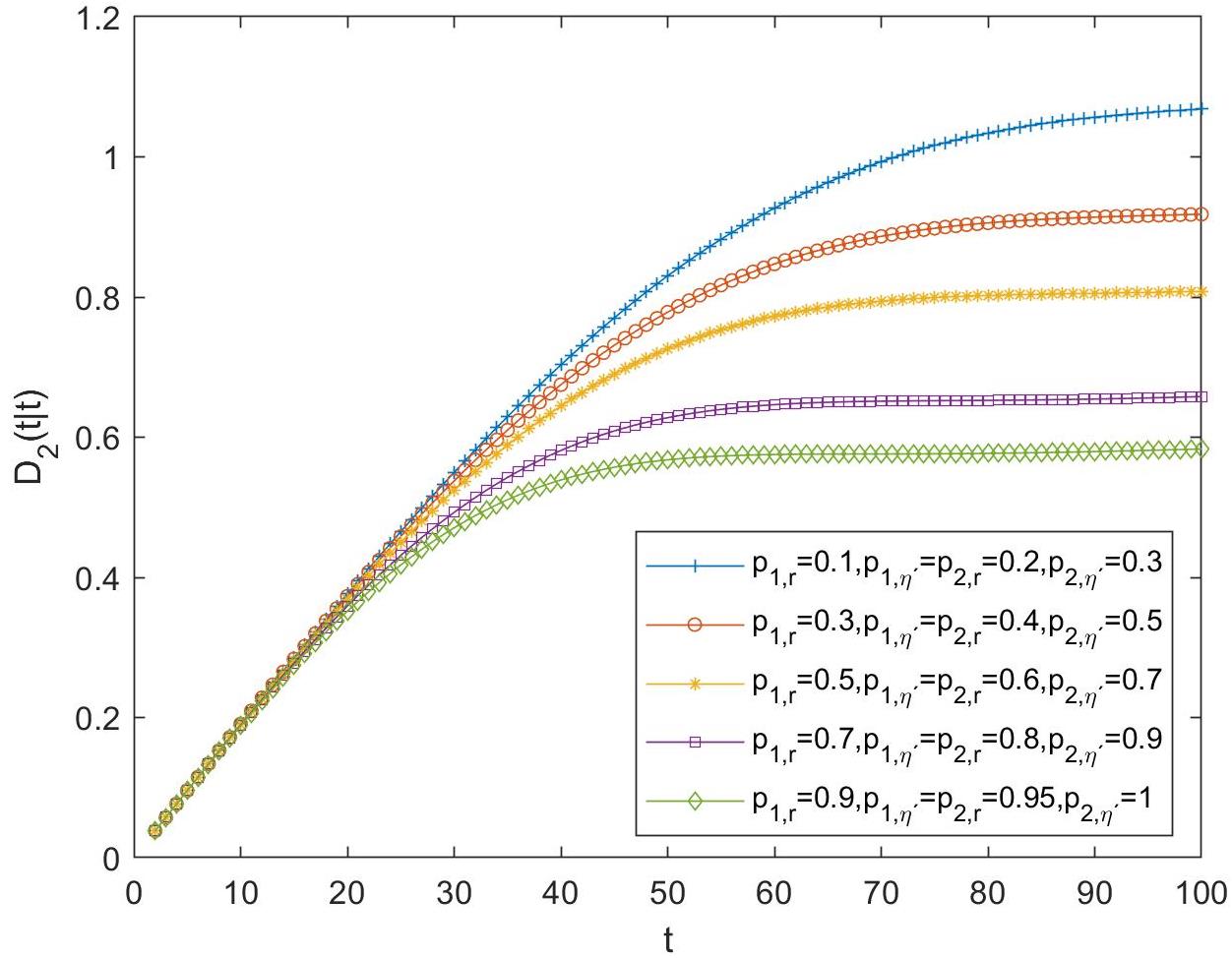}
\caption{Difference $D_2(t|t)$ between QSWL and \texorpdfstring{$\mathbb{T}_2$}-proper fusion filtering error variances for the second component of the signal for cases 16, 17, 18, 19, and 20.\label{fig8}}
\end{figure}

For all the above cases, the differences between the quaternion and tessarine filtering error variances have been calculated and displayed in Figures \ref{fig7} and \ref{fig8}, for the $\mathbb{T}_1$ and $\mathbb{T}_2$-proper scenarios and for the first and second component of the signal, respectively. Same conclusions can be derived for both figures: 1) better estimations by using the tessarine processing than from the quaternion processing and, 2) there exists a lower difference between the estimations obtained from both types of processing when the probability that the components of the available observations are updated increases.      

\section{Discussion}

The LLMS centralized fusion filtering problem  is analyzed in linear systems with multiple sensors  and multiple packet dropouts. However, unlike most of  the solutions proposed in the literature, 
a proper hypercomplex-valued signal processing has been employed with the purpose of reducing  the dimension of the problem. 
Specifically, the state-space system is defined in the tessarine domain, and it is assumed  that each  component of the measurement  output at each sensor may present a different packet dropout rate, modeled by using a Bernoulli random variable.  Moreover,   the state and the measurement noises can be correlated. Under hypotheses of $\mathbf{T}_k$-properness, our approach allows us to provide an optimal LLMS fusion filtering algorithm that reduces the computational cost of  its counterpart in the real field.  The good behavior and benefits of this algorithm has been analyzed in situations of $\mathbf{T}_1$ and  $\mathbf{T}_2$- properness, by considering different number of sensors. Moreover, a comparative study of the quaternion and tessarine approaches has been carried out, showing how the algorithm proposed behaves better than its counterpart in the quaternion domain when  $\mathbf{T}_k$-properness, $k=1,2$, conditions are satisfied. 

As a consequence, our approach based on $\mathbf{T}_k$-proper processing presents two main advantages: on the one hand, the tessarine systems offers a suitable framework to model 3D and 4D physical and experimental phenomena, and on the other hand, a considerably reduction of problem dimention is possible, when the processes involved are $\mathbf{T}_k$-proper, which implies significant computational savings in the implementation of our LLMS fusion filtering algorithm, that cannot be attained from a real formalism of the problem. 

In future research, we will approach the estimation problem in other hypercomplex algebras and under different properness conditions, by using alternative fusion architectures for the multi-sensor observations, and with varied uncertainty situations.

\authorcontributions{All authors have contributed equally to the work.     The functions mainly carried out by each specific author are detailed below.     Conceptualization, R.M.F.-A., N.L.-B., and  C.C.T.; formal analysis, R.M.F.-A., and  J.D.J.-L.; methodology, R.M.F.-A., J.D.J.-L.; investigation, R.M.F.-A., and J.D.J.-L.; visualization, R.M.F.-A., J.D.J.-L., N.L.-B., and C.C.T.; writing---original draft preparation, R.M.F.-A. and J.D.J.-L.; writing---review and editing, R.M.F.-A., J.D.J.-L., N.L.-B., and C.C.T.; funding acquisition, R.M.F.-A. and J.D.J.-L.; project administration, R.M.F.-A.; software, J.D.J.-L.; supervision, R.M.F.-A., J.D.J.-L., N.L.-B., and C.C.T.; validation, R.M.F.-A., J.D.J.-L., N.L.-B., and C.C.T. All authors have read and agreed to the published version of the~manuscript.}

\funding{This paper has been supported in part by the Project PID2021-124486NB-I00
of the Plan Estatal de I+D+i", Ministerio de Educación y Ciencia, Spain, the I+D+i project with reference number 1256911, under `Programa Operativo FEDER Andalucía 2014--2020', Junta de Andalucía, and Project EI-FQM2-2023 of `Plan de Apoyo a la Investigación 2023--2024' from the University of Jaén.}

\institutionalreview{Not applicable.}

\informedconsent{Not applicable.}

\dataavailability{Not applicable.} 


\conflictsofinterest{The authors declare no conflict of interest.} 



\abbreviations{Abbreviations}{
The following abbreviations are used in this manuscript:\\

\noindent 
\begin{tabular}{@{}ll}
LLMS & Linear least-mean-squares\\
QSL & Quaternion strictly linear\\
QSWL & Quaternion semi-widely linear\\
WL & Widely linear
\end{tabular}
}

\appendixtitles{yes} 
\appendixstart
\appendix

\section{Proof of Theorem 1}\label{Append1}

Consider the state-space system given by the equations \eqref{system3a}-\eqref{system3b},  and \eqref{y_k}. On the basis of the innovations $\boldsymbol{\varepsilon}_k(t)=\mathbf{y}_k(t)-\hat{\mathbf{y}}_k(t|t-1)$,  with $\hat{\mathbf{y}}_k(t|t-1)$ the LLMS estimator of $\mathbf{y}_k(t)$ based on the measurements $\{\mathbf{y}_k(1), \ldots, \mathbf{y}_k(t-1)\}$,  for $k=1,2$,  the LLMS filter $\hat{\bar{\mathbf{x}}}(t|t)$ of $\bar{\mathbf{x}}(t)$ can be expressed as  \cite{Navarro2020}:
\begin{equation*}\label{Local filter linear combination of the innovation}
\hat{\bar{\mathbf{x}}}(t|t)=\sum_{s=1}^t\bar{\boldsymbol{\Theta}}_k(s){\boldsymbol{\Omega}}_k^{-1}(s){\boldsymbol{\varepsilon}}_k(s),
\end{equation*}
 where $\bar{\boldsymbol{\Theta}}_k(s)=E\left[\bar{\mathbf{x}}(s){\boldsymbol{\varepsilon}}_k^{\texttt{H}}(s)\right]$, and~$ {\boldsymbol{\Omega}}_k(s)=E\left[{\boldsymbol{\varepsilon}}_k(s){\boldsymbol{\varepsilon}}_k^{\texttt{H}}(s)\right]$,
which yields the recursive equation  
\begin{equation}\label{wlfil}
\hat{\bar{\mathbf{x}}}(t|t)=\hat{\bar{\mathbf{x}}}(t|t-1)+\bar{\mathbf{L}}_k(t)\boldsymbol{\varepsilon}_k(t),
\end{equation}
with $\bar{\mathbf{L}}_k(t)=\bar{\boldsymbol{\Theta}}_k(t)\boldsymbol{\Omega}_k^{-1}(t)$.  Thus,  from the  $\mathbb{T}_k$-properness conditions, the equation  \eqref{filter} is obtained.

Moreover,  taking projections on both sides of equations \eqref{system3a} and \eqref{system3c},  we obtain that
\begin{equation}\label{wlpred}
\hat{\bar{\mathbf{x}}}(t+1|t)=\bar{\boldsymbol{\Phi}}(t)\hat{\bar{\mathbf{x}}}(t|t)+ {\bar{\mathbf{H}}_k}(t)\boldsymbol{\varepsilon}_k(t),
\end{equation}
with ${\bar{\mathbf{H}}}_k(t)=\vec{\mathbf{S}}(t)\bar{\boldsymbol{\Pi}}_k^{\vec{\boldsymbol{\gamma}}^{\texttt{H}}}(t) \boldsymbol{\Omega}_k^{-1}(t)$, and 
\begin{equation}\label{esty}
\hat{\mathbf{y}}_k(t|t-1)= \bar{\boldsymbol{\Pi}}_k^{\vec{\boldsymbol{\gamma}}}(t){\boldsymbol{\mathcal{C}}}\hat{\bar{\mathbf{x}}}(t|t-1)+ \bar{\boldsymbol{\Pi}}_k^{\left(\boldsymbol{1}-\vec{\boldsymbol{\gamma}}\right)}(t)\vec{\mathbf{y}}_k(t-1),
\end{equation}
Then, by applying  $\mathbb{T}_k$-properness conditions on \eqref{wlpred} and \eqref{esty}, and considering that $E[\bar{\mathbf{u}}(t)\boldsymbol{\varepsilon}_k^{\texttt{H}}(s)]=\vec{\mathbf{S}}(t)\bar{\boldsymbol{\Pi}}_k^{\vec{\boldsymbol{\gamma}}^{\texttt{H}}}(t)\delta_{t,s}$,  the equations \eqref{pred} and \eqref{central innovations} are directly devised.

Let $\bar{\boldsymbol{\epsilon}}(t|t-1)=\bar{\mathbf{x}}(t)-\hat{\bar{\mathbf{x}}}(t|t-1)$ be the prediction error, then
\begin{equation*}\label{wltheta}
\bar{\boldsymbol{\Theta}}_k(t)=E\left[\bar{\mathbf{x}}(t)\bar{\boldsymbol{\epsilon}}^{\texttt{H}}(t|t-1)\right]\boldsymbol{\mathcal{C}}^{\texttt{T}}\tilde{\boldsymbol{\Pi}}_k^{\vec{\boldsymbol{\gamma}}^{\texttt{H}}}(t)=
\bar{\mathbf{P}}(t|t-1) \boldsymbol{\mathcal{C}}^{\texttt{T}}\tilde{\boldsymbol{\Pi}}_k^{\vec{\boldsymbol{\gamma}}^{\texttt{H}}}(t),  \quad t\geq 2,
\end{equation*}
where $\bar{\mathbf{P}}(t|t-1)=E[\bar{\boldsymbol{\epsilon}}(t|t-1)\bar{\boldsymbol{\epsilon}}^{\texttt{H}}(t|t-1)]$. As a consequence, under $\mathbb{T}_k$-properness conditions, \eqref{Theta} is derived, where $\mathbf{P}_k(t|t-1)$ is given by the first $knR\times knR$ submatrix of $\bar{\mathbf{P}}(t|t-1)$ and also the equation \eqref{D local filter} for $\boldsymbol{\Gamma_{\bar{x}}}(t)=E[\bar{\mathbf{x}}(t)\bar{\mathbf{x}}^{\texttt{H}}(t)]$ is easily obtained from \eqref{system3a}.

In order to devise the equation \eqref{Omega central filter theorem}, we will rewrite the equation \eqref{central innovations} as follows:
\begin{equation*}\label{innoer}
\boldsymbol{\varepsilon}_k(t)=\Delta\bar{\boldsymbol{\mathcal{D}}}_k^{\vec{\boldsymbol{\gamma}}}(t)\boldsymbol{\mathcal{C}}\bar{\mathbf{x}}(t)+  \bar{\boldsymbol{\Pi}}_k^{\vec{\boldsymbol{\gamma}}}(t)\boldsymbol{\mathcal{C}}\bar{\boldsymbol{\epsilon}}(t|t-1) +\bar{\boldsymbol{\mathcal{D}}}_k^{\vec{\boldsymbol{\gamma}}}(t)\vec{\mathbf{v}}(t)+ \Delta \bar{\boldsymbol{\mathcal{D}}}_k^{(\boldsymbol{1}-\vec{\boldsymbol{\gamma}})}(t)\vec{\mathbf{y}}(t-1).
\end{equation*}
Then, considering that $\vec{\mathbf{v}}(t)$ is orthogonal to $\bar{\mathbf{x}}(t)$,  $\vec{\mathbf{y}}(t-1)$, and $\bar{\epsilon}(t|t-1)$, and $E[\Delta\bar{\boldsymbol{\mathcal{D}}}_k^{\vec{\boldsymbol{\gamma}}}(t)]=0$, 
we have that 
\begin{equation}\label{Omega_proof}
\begin{aligned}
\boldsymbol{\Omega}_k(t)=E\left[\Delta\bar{\boldsymbol{\mathcal{D}}}_k^{\vec{\boldsymbol{\gamma}}}(t)\boldsymbol{\mathcal{C}}\bar{\mathbf{x}}(t)\bar{\mathbf{x}}^{\texttt{H}}(t)\boldsymbol{\mathcal{C}}^{\texttt{T}} \Delta\bar{\boldsymbol{\mathcal{D}}}_k^{\vec{\boldsymbol{\gamma}}^{\texttt{H}}}(t)\right]+
E\left[\Delta\bar{\boldsymbol{\mathcal{D}}}_k^{\vec{\boldsymbol{\gamma}}}(t)\boldsymbol{\mathcal{C}}\bar{\mathbf{x}}(t)\vec{\mathbf{y}}^{\texttt{H}}(t-1)\Delta\bar{\boldsymbol{\mathcal{D}}}_k^{(\boldsymbol{1}-\vec{\boldsymbol{\gamma}})^{\texttt{H}}}(t)\right]\\ +E\left[\Delta\bar{\boldsymbol{\mathcal{D}}}_k^{(\boldsymbol{1}-\vec{\boldsymbol{\gamma}})}(t)\vec{\mathbf{y}}(t-1)\bar{\mathbf{x}}^{\texttt{H}}(t)\boldsymbol{\mathcal{C}}^{\texttt{T}} \Delta\bar{\boldsymbol{\mathcal{D}}}_k^{\vec{\boldsymbol{\gamma}}^{\texttt{H}}}(t)\right]
+E\left[\bar{\boldsymbol{\mathcal{D}}}_k^{\vec{\boldsymbol{\gamma}}}(t)\vec{\mathbf{v}}(t)\vec{\mathbf{v}}^{\texttt{H}}(t)\bar{\boldsymbol{\mathcal{D}}}_k^{\vec{\boldsymbol{\gamma}}^{\texttt{H}}}(t)\right]
\\ + E\left[\Delta\bar{\boldsymbol{\mathcal{D}}}_k^{(\boldsymbol{1}-\vec{\boldsymbol{\gamma}})}(t)\vec{\mathbf{y}}(t-1)\vec{\mathbf{y}}^{\texttt{H}}(t-1)\Delta\bar{\boldsymbol{\mathcal{D}}}_k^{(\boldsymbol{1}-\vec{\boldsymbol{\gamma}})^{\texttt{H}}}(t)\right] + \bar{\boldsymbol{\Pi}}_k^{\vec{\boldsymbol{\gamma}}}(t)\boldsymbol{\mathcal{C}}\bar{\mathbf{P}}(t|t-1)\boldsymbol{\mathcal{C}}^{\texttt{T}}\bar{\boldsymbol{\Pi}}_k^{\vec{\boldsymbol{\gamma}}^{\texttt{H}}}(t).
\end{aligned}
\end{equation}

The equation \eqref{Omega central filter theorem} follows from \eqref{Omega_proof}, by using Hadamard product properties and taking into account that
$$\boldsymbol{\Gamma}_{\bar{\mathbf{x}}\vec{\mathbf{y}}}(t,t-1)=\bar{\boldsymbol{\Phi}}(t-1){\boldsymbol{\Gamma}}_{\bar{\mathbf{x}}\vec{\mathbf{y}}}(t-1) + \vec{\mathbf{S}}(t-1) \bar{\boldsymbol{\Pi}}^{\vec{\gamma}}(t-1). $$

Finally, from \eqref{wlfil}, it is clear that the pseudo covariance matrix  $\bar{\mathbf{P}}(t|t)=E[\bar{\boldsymbol{\epsilon}}(t|t)\bar{\boldsymbol{\epsilon}}^{\texttt{H}}(t|t)]$, of the filtering errors $\bar{\boldsymbol{\epsilon}}(t|t)=\bar{\mathbf{x}}(t)-\hat{\bar{\mathbf{x}}}(t|t)$, can be computed of the form 
$$\bar{\mathbf{P}}(t|t)=\bar{\mathbf{P}}(t|t-1) - \tilde{\boldsymbol{\Theta}}_k(t)\boldsymbol{\Omega}_k^{-1}(t)\tilde{\boldsymbol{\Theta}}_k^{\texttt{H}}(t),$$
and consequently, using the $\mathbb{T}_k$ properness conditions, \eqref{centralerrorfil} holds.
Furthermore, since 
\begin{equation}\label{errorpred}
\bar{\boldsymbol{\epsilon}}(t+1|t)=\bar{\mathbf{x}}(t+1)-\hat{\bar{\mathbf{x}}}(t+1|t)=\bar{\boldsymbol{\Phi}}(t) \bar{\boldsymbol{\epsilon}}(t|t)+\bar{\mathbf{u}}(t)-\tilde{\mathbf{H}}_k(t)\boldsymbol{\varepsilon}_k(t),
\end{equation}
$E[\boldsymbol{\epsilon}(t|t) \boldsymbol{\varepsilon}_k(t)]=0$, and $E[\bar{\mathbf{u}}(t) \boldsymbol\bar{\boldsymbol{\epsilon}}^{\texttt{H}}(t|t)]=-\vec{\mathbf{S}}(t)\bar{\boldsymbol{\Pi}}_k^{\vec{\boldsymbol{\gamma}}}(t)\bar{\mathbf{L}}_k^{\texttt{H}}(t) $, we obtain that
\begin{multline*}\label{WLprederror}
\bar{\mathbf{P}}(t+1|t)=\bar{\boldsymbol{\Phi}}(t)\bar{\mathbf{P}}(t|t)\bar{\boldsymbol{\Phi}}^{\texttt{H}}(t)- \bar{\mathbf{H}}_k(t)\bar{\boldsymbol{\Theta}}_k^{\texttt{H}}(t)\bar{\boldsymbol{\Phi}}^{\texttt{H}}(t)\\- \bar{\boldsymbol{\Phi}}(t) \bar{\boldsymbol{\Theta}}_k(t) \bar{\mathbf{H}}_k^{\texttt{H}}(t)- \bar{\mathbf{H}}_k(t) \boldsymbol{\Omega}_k(t) \bar{\mathbf{H}}_k^{\texttt{H}}(t)+\bar{\mathbf{Q}}(t),
\end{multline*}
and hence, from $\mathbb{T}_k$ properness conditions, \eqref{centralerrorpred} follows.


\begin{adjustwidth}{-\extralength}{0cm}

\reftitle{References}\label{references}

\end{adjustwidth}

\end{document}